\documentclass[12pt,reqno]{amsart}
\usepackage[english]{babel}
\usepackage[colorlinks,linkcolor=blue]{hyperref}
\usepackage[nameinlink]{cleveref}
\usepackage{amsmath,amssymb,bbm}
\usepackage{mathtools}
\usepackage{graphicx, color}
\usepackage[dvipsnames]{xcolor} 
\usepackage{enumitem}  
\usepackage{pgf} 
\usepackage{geometry}
\makeatletter
\@namedef{subjclassname@2020}{\textup{2020} Mathematics Subject Classification}
\makeatother

\theoremstyle{plain}
\begingroup
\theoremstyle{plain}
\newtheorem{theorem}{Theorem}[section]

\newtheorem{proposition}[theorem]{Proposition}
\newtheorem{lemma}[theorem]{Lemma}

\theoremstyle{remark}
\newtheorem{remark}[theorem]{Remark}
\endgroup

\theoremstyle{definition}
\theoremstyle{remark}

\numberwithin{equation}{section}
 
\newcommand{\Rz}{{\mathbb R}}

\newcommand{\haz}{\widehat}

\newcommand{\epsi}{\varepsilon}
\newcommand{\ove}{\overline}

\renewcommand{\d}{{\rm d}}

\newcommand{\mres}{\mathbin{\vrule height 1.6ex depth 0pt width
    0.13ex\vrule height 0.13ex depth 0pt width 1.3ex}}
\newcommand{\zz}{z}

\title[Free boundary growth problem]{A free boundary problem in accretive growth}

\author[U. Stefanelli]{Ulisse Stefanelli} 
	\address[Ulisse Stefanelli]{University of
		Vienna, Faculty of Mathematics,
                Oskar-Morgenstern-Platz 1, A-1090 Vienna, Austria.}
	\email{ulisse.stefanelli@univie.ac.at}
	\urladdr{http://www.mat.univie.ac.at/$\sim$stefanelli}


\begin{document}

\begin{abstract} We  study  a free boundary
  problem inspired by the  modeling  of accretive growth. The growth
  process is formulated through a level-set approach, leading to a
  boundary-value problem for a Hamilton--Jacobi equation within a
  prescribed constraining set. Existence, variational representability, and regularity of
  solutions to the growth subproblem are investigated. The full system
  arises from coupling the growth dynamics with an elliptic equation
  for the activation field. Existence of solutions to the fully
  coupled free boundary  problem is obtained via
  a fixed-point argument. 

\end{abstract}

\maketitle

\section{Introduction}

We are interested in the free boundary problem 
\begin{align} 
  &H(x,Ku(x,v(x)),\nabla v)=0\quad \text{in} \ \Omega\setminus
    \ove{V_0},\label{eq:02}\\
  &-\Delta u(x,t) =1\quad \text{for} \ x\in  \{v <t\}, \ t>0.\label{eq:01}
\end{align}
Here, $\Omega\subset \Rz^n$ is  possibly 
unbounded open set,  and  $V_0 \subset \Omega$.
The unknowns are the two scalar fields $v: \ove \Omega  \to
[0,\infty)$ and $u : Q_v \to \Rz$, where the latter  is defined
on the a priori unknown space-time noncylindrical
domain $Q_v:=\{(x,t)\in  \Omega \times (0,\infty) \ | \ v(x)<t\}$. The
Hamiltonian $H:\ove \Omega \times \Rz \times \Rz^n \to \Rz$ is continuous,
with $p \mapsto H(x,u,p)$ having compact, convex, and nondegenerate
 zero sublevel, 
and $K$ is a space-time convolution operator.

The main result of the paper, Theorem \ref{thm:main}, states that
there exists a solution  $(u,v)$  to \eqref{eq:02}--\eqref{eq:01}
under mixed Dirichlet-Neumann conditions.  The stationary Hamilton--Jacobi equation \eqref{eq:02} is solved in
the viscosity sense, while the elliptic equation \eqref{eq:01} is
treated in the weak sense.  Along the way, we also discuss existence, variational
representability, and regularity of sublevels for viscosity solutions
of \eqref{eq:02} for a given $u$, see 
Proposition \ref{prop:representation} and Proposition~\ref{prop:john}.

\subsection{Accretive growth}
Problem \eqref{eq:02}--\eqref{eq:01} is inspired by the modelization
of an {\it accretive-growth}
process in a confined environment. 
Accretive growth occurs across a variety of
different biological systems, as well as in many natural and
technological settings.  The formation of horns, teeth, and seashells
\cite{skalak,thompson}, coral reefs \cite{Kaandorp}, bacterial
colonies \cite{Hoppensteadt},  trees
\cite{fournier-et-al}, and cell motility due to actin growth
\cite{hodge} are biological examples of accreting systems, among many others
\cite{goriely}. In geophysics, sedimentation \cite{Bustos} and glacier
evolution \cite{Hutter} are also accretive processes, as well as planet
formation \cite{Safronov}. 
 Furthermore,  accretive growth is a key aspect in metal solidification
\cite{schwerdtfeger-et-al}, crystal growth \cite{langer},  
additive manufacturing~\cite{horn}, layering, coating, and masonry,
just to mention a few.

The growing  system  is assumed to be constrained to a given smooth set
$\ove \Omega
\subset \Rz^n$, which may be chosen to be unbounded. 
We follow the classical {\it level-set
approach} \cite{Evans-Spruck,Souganidis}. Starting from the initial
shape $V_0
\subset \Omega$, the growing system is described by the set-valued function $t \mapsto V(t)
\subset \Rz^n$, whose values are the sublevels of the
unknown  function  $v:\ove \Omega \to [0,\infty)$, namely, 
$$V(t) := \{x\in \Omega\ | \ v (x)<t\}\quad \text{for}
\ t>0.$$
The value $v(x)$ models the {\it time} at which the point $x \in
 \Omega $ is added to the growing system. Correspondingly,
$v$ is called  {\it time-of-attachment} function in the
following. 

Growth is always influenced by its environment and aspects such as
nutrient or oxygen concentration, pressure, temperature, or stress
\cite{goriely}, often have to be taken into account. We model this fact by
considering the simplest situation of a scalar  {\it activation field} $u:Q_v \to \Rz$
driving the growth. The dynamics of the activation
field $u$ is influenced by the actual shape of the growing
system. Since growth is often very slow with respect to other
processes, we assume that $u$ evolves quasistatically. We model  a
basic quasistatic
diffusion process by
the linear elliptic equation \eqref{eq:01} on the growing domain
$V(t)$ for all $t>0$. This corresponds to the case of $u$ being a
concentration or temperature in a homogeneous medium and could  of course be  generalized. In particular, the choice of the right-hand side equal
to $1$ is made solely for simplicity.
Note that the resulting full model 
\eqref{eq:02}--\eqref{eq:01}  is of free boundary type.

Starting from the initial shape $V_0$, the growing system
$V(t)$   is  assumed to evolve by layering material at its
boundary at a
velocity depending on the position $x \in \partial V(t)$, on the outward-pointing normal
direction $n(x)$  at  $x\in  \partial V(t)$, and on the activation field $u$.  The evolution in time
of a 
point $x(t) \in \partial V(t)$ is modeled  by the ODE
$$ x'(t)=f\big(x(t),u(x(t),t),n(x(t))\big),$$ where $f: \Rz^n \times \Rz \times
\Rz^n\to \Rz^n$ is a given velocity field and the prime indicates time
 differentiation.    
Assuming smoothness of $\partial V(t)$ (an assumption
which we will soon be forced to  abandon, however),
differentiating $t=v(x(t))$ with respect to time, and setting 
$n(x)=\nabla v(x)/|\nabla v(x)|$, the function $v$ is
then asked to satisfy $H(x,u(x,v(x)),\nabla v(x))=0$ for the
Hamiltonian
$H(x,u,p)=  p\cdot f(x,u,p/|p|)-1$.  A reference situation is that
of {\it normal accretion}, where $f(x,u,n)= \gamma(x,u,n) n$ for some
scalar accretion rate $\gamma: \Rz^n \times \Rz \times \Rz^n \to 
(0,\infty)$. In this case, the Hamiltonian reads $H(x,u,p) =  
\gamma(x,u,p/|p|)|p|-1$. We make no assumption on the specific
form of the Hamiltonian in the following, but simply ask $H$ to be
continuous and the sets $\{p\in \Rz^n\,|\, H(x,u,p)\leq 0\}$ to be
convex, bounded, and nondegenerate.  
 
A crucial trait of the model is the occurrence of the term
$u(x,v(x))$ in the
Hamiltonian $H$ in \eqref{eq:02}. This reflects the assumption that the growth
process is driven by the value
of the activating field {\it at the growing front}. On the other hand,
treating the
composition $x \mapsto u(x,v(x))$ of the two unknowns $u $ and $v$ into the
Hamiltonian $H$ poses severe
mathematical challenges. In fact, $V(t)$ may develop
singularities, see Figure \ref{fig:1}, and can change topology.
As a consequence, the activation field
$u$ cannot be assumed to be regular in time, which prevents one from
taking its trace on the a priori unknown and nonsmooth
space-time manifold $\{(x,t)\in  \ove \Omega \times
(0,\infty)\ | \ t= v(x) \}$. We are hence forced to regularize
the problem by replacing $u$   in the Hamilton--Jacobi equation
\eqref{eq:02}  by a space-time convolution $K u$. This provides the
additional regularity needed to take the trace of the field $u$ on
$\{(x,t)\in \ove \Omega  \times
(0,\infty)\ | \ t= v(x) \}$. 
From the modeling viewpoint,
this relaxation implies that the growth is driven by a local mean of
the values of
the activating field $u$ at the accreting front, instead of its
pointwise ones. This modification has a limited impact on the
modeling, as we are entitled to choose the support  of the
regularizing kernels
to be arbitrarily small, so  as  to localize the
regularization.

\subsection{Boundary conditions}
\begin{figure}
  \centering
  \pgfdeclareimage[width=95mm]{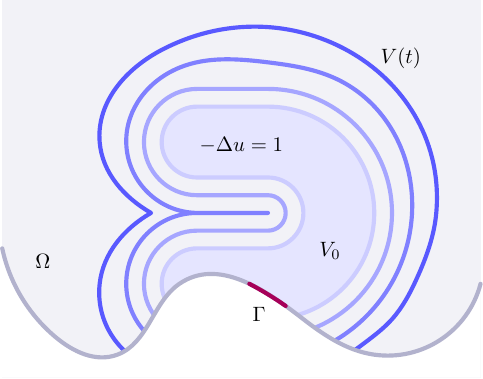}{figure}
  \pgfuseimage{figure}
  \caption{Setting and notation.}
  \label{fig:1}
\end{figure}

The PDE system
\eqref{eq:02}--\eqref{eq:01} calls for specifying boundary
conditions. As $V_0 \subset \Omega$ models the {\it initial} shape
and $v$ is the time-of-attachment function,
the Hamilton--Jacobi equation \eqref{eq:02} is naturally
complemented by  the condition  
\begin{equation}
   v =0 \quad \text{on} \ V_0.\label{eq:044}
 \end{equation}
 At the boundary $\partial \Omega$ of the constraining set we can
 consider the boundary condition  
 \begin{equation}
   \label{eq:055}
   H(x,Ku(x,v(x)),\nabla v) \geq 0 \quad \text{for} \ x \in
   \partial \Omega \setminus \ove{V_0}.
 \end{equation}
 In the setting of viscosity solutions to Hamilton--Jacobi constrained
 problems, this boundary condition is classical and qualifies $v$
 as {\it constrained} viscosity solution, see  
 \cite{CapuzzoDolcetta-Lions,Soner1,Soner2}.

The activation field $u$ is subject to mixed homogeneous
 Dirichlet--Neumann conditions. 
To pose boundary conditions on $u$ is subtle, as
$\partial V(t)$ need not be smooth. The Dirichlet condition is hence solely
imposed on $\Gamma \subset \partial  V_0 \cap \partial \Omega$,
which is assumed to be nonempty, open (with respect to the topology of $\partial \Omega$), and Lipschitz. As $ V(t)$ grows
with time, we have that
$\Gamma \subset \partial V(t)$ for all $t>0$ and we can prescribe
\begin{align}
  &u =0 \quad \text{on} \ \Gamma \times (0,\infty).\label{eq:033}
\end{align} 
 The homogeneous Neumann condition is then imposed on the complement
 $\partial V(t)\setminus \ove \Gamma$, namely,
 \begin{align}
  &\partial u/\partial n =0 \quad \text{on} \ \partial V(t)\setminus
    \ove \Gamma, \   \ t>0.\label{eq:066}
\end{align}
Note nonetheless that condition \eqref{eq:066} will be imposed just
weakly: We address a variational formulation of equation
\eqref{eq:01}, see \eqref{eq:0111} below, ensuring the validity of the
Neumann boundary condition where $\partial
V(t)$ is smooth enough.  

 Let  us remark that the  Neumann boundary condition
\eqref{eq:066} is also imposed on the portion  $A(t):=(\partial V(t)\setminus \ove \Gamma)
\cap \partial \Omega$, which  under our assumptions will be  nonempty, see Figure \ref{fig:1}. From the
modeling viewpoint, this is  somewhat  still unsatisfactory, as one would
prefer to have the Dirichlet condition there. This would however
require $A(t)$ to be sufficiently regular in order for a trace to be
well defined and such regularity of $A(t)$ is presently not known.

\subsection{Results} Our main result, Theorem
\ref{thm:main}, establishes the existence of suitably weak solutions
  $(u,v)$   to \eqref{eq:02}--\eqref{eq:066}. More precisely, the
Hamilton--Jacobi equation \eqref{eq:02}, subject to conditions
\eqref{eq:044}--\eqref{eq:055}, is solved in the viscosity
sense, whereas the elliptic equation \eqref{eq:01} on the unknown
domain $V(t)$, under conditions \eqref{eq:033}--\eqref{eq:066}, is
treated  variationally. 

Existence of solutions is proved by an iterative argument, where we
 alternately solve the Hamilton--Jacobi subproblem \eqref{eq:02}, \eqref{eq:044}--\eqref{eq:055} for a given 
$u$, and the elliptic subproblem \eqref{eq:01},
\eqref{eq:033}--\eqref{eq:066} for a given $v$.

In addition to accretive growth, the Hamilton--Jacobi subproblem \eqref{eq:02},
\eqref{eq:044}--\eqref{eq:055} (for $u$ given) arises in different contexts, from front
propagation and geometric optics \cite{Barles2,Evans-Spruck,Souganidis}, to
optimal control \cite{BardiCapuzzoDolcetta}, differential games
\cite{Cardaliaguet,Gomes}, and image reconstruction \cite{Aubert}.
As a consequence, it has attracted
early attention, see for instance the
classical references
\cite{CapuzzoDolcetta-Lions,userguide,Lions}. Different results are
available, corresponding to different assumptions on the ingredients of
the problem. Still, we have not been able to find a result which exactly covers  the
specific framework of \eqref{eq:H01}--\eqref{eq:H04} below.
For this reason, we present complete arguments,  not  necessarily
claiming originality.
Existence and variational representation of nonnegative viscosity
solutions to \eqref{eq:02},
\eqref{eq:044}--\eqref{eq:055} (for $u$ given) are proved in 
Proposition~\ref{prop:representation}.

The solution to the elliptic subproblem \eqref{eq:01},
\eqref{eq:033}--\eqref{eq:066} for a given $v$ relies on the
possibility of establishing a Poincar\'e  inequality  on
$V(t)$. In turn, this requires proving that $V(t)$ maintains some minimal
level of regularity. In Proposition~\ref{prop:john}, we
show that, starting from some suitably regular $V_0$, the
growing set $V(t)$ is a {\it John domain} for all times (see Section
\ref{sec:main}). This extends to the constrained case $\Omega\not
=\Rz^n$ the arguments in  \cite{john}. Some   technical modifications
are required for this extension, since the boundary 
$\partial \Omega$ of the constraining domain must also be taken into
account. Note that John regularity is sharp in this setting \cite{john}.

 Finally,   the convergence of the iterative scheme depends on the time
regularity of the set-valued maps $t \mapsto V(t)$ for different
$v$. These are Lipschitz continuous with respect to the
Hausdorff distance, as well as H\"older continuous with respect to Lebesgue
measure of the differences. Crucially, such
continuities are uniform with respect to the iteration, 
allowing for the use of compactness arguments.   

Let us conclude this section by noting that uniqueness of solutions to the coupled system
\eqref{eq:02}--\eqref{eq:066} remains open. Even under the boundary condition
\eqref{eq:055}, classical uniqueness arguments
\cite{Barles,CapuzzoDolcetta-Lions,Lions} seem to be 
inapplicable. In fact, without further assumptions uniqueness seems to fail even in
the simplest decoupled setting of  a nonconvex  $H$ not depending on $u$, see
Remark \ref{rem:uniqueness} below.   Uniqueness holds however for $H$
convex and Lipschitz under condition \eqref{eq:055}. 

Before moving on, let us mention that the  setting of
\eqref{eq:02}--\eqref{eq:066} is not meant to target maximal 
generality  
but rather aimed at introducing a class of problems and its
challenges. In particular, the  results of this paper can be
generalized to more complex  settings.  A first generalization consists in considering the Hamilton--Jacobi equation
$$H(x,v(x),Ku(x,v(x)), \nabla v)=0\quad x \in 
\Omega\setminus \ove{V_0}$$
where an additional dependence on $v$ is included. This
extension would allow to cover the case of a time-dependent growth
process at the cost of a somewhat more intricate analysis. A
second possible direction is to replace   equation
\eqref{eq:01} with other stationary PDEs, provided the driving operator retains suitable compactness and continuity properties. In particular, \eqref{eq:01} can be replaced by a semilinear elliptic system with sufficiently controlled nonlinearities.

\subsection{Connections with previous studies}
The theory of the Hamilton--Jacobi subproblem \eqref{eq:02} with condition
\eqref{eq:044}
for $u$ given (and suitably smooth)  is well developed. We limit
ourselves  to  
mentioning the classical references
\cite{BardiCapuzzoDolcetta,Barles,userguide,Lions} for the general
frame of viscosity solutions. In the context of optimal control
with state constraints, condition \eqref{eq:055} corresponds to the
notion of {\it constrained} viscosity solution, going back to
\cite{CapuzzoDolcetta-Lions,Soner1,Soner2}, see also
\cite{Altarovici,Frankowska,Ishii-Koike,Kim,Loreti} among many
others. A class of biological invasion models based on
Hamilton--Jacobi dynamics with obstacles has been introduced in \cite{Bressan0}
and studied in \cite{Bressan,DeLellis}, see also
\cite{Bressan2,Bressan3} and the reference therein.
Our results rely fundamentally on the metric structure
underlying the Hamilton--Jacobi subproblem, which is already pointed out
in  \cite{Courant-Hilbert2,Lions}, see also
\cite{Mennucci,Siconolfi}. A wide range of results on the
Hamilton--Jacobi subproblem is already available.

When considering the full free boundary problem, the closest
paper to the present one  is~\cite{morpho3}. This focuses on
\eqref{eq:02}--\eqref{eq:01}  for   the
specific {\it generalized eikonal} case $H(x,u,p)= \gamma(u)|p|-1$ and
for \eqref{eq:01} replaced by the linearized elastic equilibrium
system. The setting in \cite{morpho3} is however much simplified by
the absence of a  constraining   set, i.e., $\Omega=\Rz^n$, and by the
choice of the docking condition $u=0$ on $\omega\times (0,\infty)$ with
$\omega \subset V_0$ instead
of the Dirichlet condition \eqref{eq:033}.  The same docking condition could
be considered  here,  as well, and would lead to a simpler theory. Note nonetheless that the
applicative relevance of  this  condition is limited to $n=2$,
as  one cannot practically realize such a docking  for  $n\geq 3$.

A related set of questions on tumor-growth modeling and
Hamilton--Jacobi dynamics is addressed by {\sc Collins,  Jacobs, \&
  Kim} in \cite{Collins}, although in the context of a different
model. There, the focus lies on the regularity of the free
boundary of a tumor model driven by pressure and by nutrient
concentration. The  central object  of the analysis is again the
function $v$, which is proved to be H\"older continuous with
regular sublevels, except possibly on a subset of Hausdorff dimension
$n-\mu$ for some $\mu \in(0,1)$. Although the underlying dynamics essentially differ, these results bear a
close analogy to Proposition~\ref{prop:john} below. A previous
analysis of the same model in specific parameter regimes is in \cite{Jacobs,Kim2}.

More generally, one can refer to the extensive literature on  
Hele--Shaw/Muskat/Stefan dynamics, see, e.g.,
\cite{Caffarelli,Visintin}.
In these models, however, growth is driven by a fundamentally
different underlying mechanism. In the Hele--Shaw setting, the driving
field satisfies a Dirichlet boundary condition, and the domain evolves with a normal velocity proportional to the gradient of this field at the free boundary. In contrast, in our framework the activating field $u$
 is subject to a homogeneous Neumann condition, while the velocity of
the free boundary depends on the values of $u$ at the growing front.
 
In a slightly different framework, it is also worth mentioning the results by {\sc Su
  \& Burger} \cite{Su} and {\sc Cardaliaguet, Ley, \& Monteillet}
\cite{Cardaliaguet2} concerning a growth model for polymers. In these papers, no
constraining set is imposed, i.e., $\Omega=\Rz^n$, and the
activating field is assumed to fulfill a parabolic equation driven by
the unknown evolving set $V(t)$. More precisely, in \cite{Su} the equation $u_t-\Delta u =-\partial_t h(v-t)$ is
considered, where $h$ is the Heaviside graph, while in \cite{Cardaliaguet2} the authors
study 
$u_t-\Delta u = g(u){\mathcal
  H}^{n-1}\mres \partial V(t)$  where the latter stands for the
$(n-1)$-dimensional Hausdorff measure restricted to $\partial V(t)$. 
The increased complexity of the nonlinear dynamics of the activating field is compensated by a
comparatively simpler treatment of the Hamilton--Jacobi subproblem,
where one assumes $H(x,u,p)=\gamma(x,u)|p|-1$. Since $u$
is defined everywhere in $\Rz^n\times (0,+\infty)$, problem
\eqref{eq:02} can be conveniently
reformulated in terms of  $w=w(x,t)$ solving
  $-w_t+\gamma(x,u(x,t))|\nabla w|=0$, so that no explicit composition of the
variables $u$ and $v$, as in  \eqref{eq:02}, is needed. Notably,
the analyses in \cite{Cardaliaguet2,Su} rely on the possibility of
proving that the front $\partial V(t)$ is Lipschitz continuous, a
regularity that cannot be expected in our setting \cite{john}.

\subsection{Structure of the paper}  In Section~\ref{sec:main},
we  introduce notation, specify assumptions, and state the main
result, i.e.,  Theorem~\ref{thm:main}. Section \ref{sec:setting} is devoted to the
study of the Hamilton--Jacobi subproblem.  Finally,  the proof of
Theorem~\ref{thm:main} is presented in Section~\ref{sec:elastic}.

\section{Setting and main result}\label{sec:main}

In this section, we fix the assumptions on the ingredients of the
problem, present the notion of solution, and state our main result,
Theorem~\ref{thm:main}.

In the following, we  set   $\Rz_+=(0,\infty)$  and  denote by $B_r(x)$ the
open ball with center $x \in \Rz^n$ and radius $r>0$, simply 
writing 
$B_r$ for $B_r(0)$.  We use the same symbol $|\cdot|$ for the Euclidean norm in $\Rz^m$ for $m\geq 1$ and for the Lebesgue measure of a
measurable set. We let $\ove A$,
$A^\circ$, $A^c$, and $1_A$ be the closure, the interior, the complement
in $\Rz^n$, and the characteristic function of a set $A \subset \Rz^n$,
respectively. Specifically, $1_A(x) = 1$ if $x \in A$ and $1_A(x)=0$
if $x \in A^c$.   
We indicate by
$d(x,A)=\inf_{a\in A}|x-a|$ the distance  from  a nonempty set $A\subset
\Rz^n$. For any two  
sets $A,\,B \subset \Rz^n$, the symbol $A+B$ stands for the usual sum.
In case $A$ and $B$ are nonempty, $d_H(A,B)$ denotes their
{\it Hausdorff distance},  namely $d_H(A,B):= \max\{e(A,B),e(B,A)\}$
where $e(A,B):= \sup_{a\in A}d(a,B)$ is the {\it excess} of $A$ with
respect to 
$B$.  If $A \subset B$, the notation ${\rm int}_B A$ indicates the
interior of $A$ with respect to $B$. 

\subsection{Assumptions} The growth process is constrained to a given
\begin{equation}
  \label{eq:H05} \Omega \subset \Rz^n\quad \text{nonempty, open, 
    connected, and  uniformly  $C^{1,1}$}.
  \end{equation}
This regularity implies that all points $x\in \partial \Omega$ can be touched by
an internal and  by  an external ball with  sufficiently small
 radius $\ove \epsi>0$.  Namely, \eqref{eq:H05} gives 
\begin{align*}
  & \exists \ove \epsi > 0 \ \forall x \in \partial \Omega \ \exists
    x_i ,\, \exists x_e \in
 \Rz^n:\\
&\quad B_{\ove \epsi}(x_i) \subset \Omega, \quad B_{\ove \epsi}(x_e) \subset  ({\ove \Omega})^c , \quad \text{and} \ \ |x-x_i|=|x-x_e|=\ove
                            \epsi,
\end{align*}
see \cite[Cor.~2]{Lewicka}.  In particular,  the projection $\pi:\partial \Omega +
B_{\ove \epsi} \to \partial \Omega$ given by
$$\pi(x):={\rm arg\, \min}_{ y \in \partial \Omega}|x-y|,$$
as well as the normal
outward-pointing vector field $\nu: \partial \Omega \to {\mathbb
  S}^{n-1}=\{ \xi  \in \Rz^n |\, | \xi |=1\}$
are well-defined and Lipschitz continuous  on $\partial \Omega$.  In the  proof of
Proposition~\ref{prop:john},   the
Lipschitz continuous composition $\eta:=-\nu\circ \pi:  (\partial \Omega +
B_{\ove \epsi/2} )  \cap \ove \Omega \to {\mathbb S}^{n-1}$ will play a
role.

As $\Omega$ is path-connected and $C^{1,1}$, its closure $\ove \Omega$ is
path-connected, as well. In particular, $\ove \Omega$ is  a geodesic space. For later use, we indicate by $d_{\ove \Omega}$ the  geodesic distance
on $\ove \Omega$, namely,
$$d_{\ove \Omega}(x,y) = \min\left\{\int_0^1 |\gamma'(s)|\, \d s \ \Big| \ \gamma \in
W^{1,\infty}(0,1;\ove \Omega), \ \gamma(0)=y,\ \gamma(1)=x\right\}\quad \forall
x,\, y \in \ove \Omega.$$
Given $x,\,y\in \ove \Omega$, we call {\it geodesic} from $y$ to $x$
 any   curve
realizing the minimum above, recalling that geodesics are not necessarily
unique. We additionally ask that the Euclidean and the geodesic
 distances  are comparable in $\ove \Omega$, namely,
\begin{equation} 
  \exists L>0: \quad d_{\ove \Omega}(x,y) \leq L|x-y| \quad \forall x,\,
    y \in \ove \Omega.\label{eq:E2} 
  \end{equation}
  Note that this does not follow from  uniform  $C^{1,1}$ regularity, as $\Omega$ may
  be unbounded.

Let us recall that a bounded set $U \subset
\Rz^{ n}$ is said to be  a  
{\it John domain} (or {\it John regular}) if there exists a point $x_0\in U$ ({\it John
  center}) and a constant $\kappa \in (0,1]$ ({\it John constant})
such that for all points $x\in U$ one can find an
arc-length parametrized
curve $\gamma: [0,\ell_\gamma] \to U$ with $\gamma (0)=x$,
  $\gamma(\ell_\gamma) = x_0$, and $ d  (\gamma (s), \partial U) \geq
  \kappa s$ for all $s\in [0,\ell_\gamma]$, see \cite{Martio}.  In
  other words, the {\it twisted cone} $\cup_{s\in
    [0,\ell_\gamma]} B_{\kappa s}(\gamma(s))$ is asked to be contained in $U$.  
  Note that John domains are connected.
  
   We pose the following assumption on  the initial set $V_0$ 
 \begin{align}
  & V_0 \subset \Omega \quad \text{is nonempty, open, and John
    regular with respect}\nonumber\\
   &\quad \text{to $x_0 \in V_0$ with John constant $\kappa_0 \in
   (0,1]$}. \label{eq:H00u}
 \end{align}
Note that, given \eqref{eq:H00u}, $\ove{V_0}$ is  John regular  with
respect to $x_0  $ with John constant $\kappa_0$, as well. 
In the following, we indicate by $ {\mathcal C}_x $ the set of Lipschitz curves
  connecting $\ove{V_0}$ to $x \in \Rz^n$, namely,
$$ {\mathcal C}_x =\{\gamma \in W^{1,\infty}(0,1;\Rz^n)\ | \ \gamma(0)\in \ove{V_0},\
\gamma(1)=x\} \quad\forall x\in \Rz^n.$$

In order to specify boundary conditions for
\eqref{eq:01} we assume the following
\begin{equation}
  \label{eq:H07}
  \Gamma \subset \partial V_0\cap \partial \Omega  \ 
  \text{nonempty, open in $\partial \Omega$, with
    $ \inf_{g\in \Gamma}  d(g,\partial V_0\setminus \partial \Omega)>0$}.
\end{equation}
As $V(t) \subset \Omega$ grows with $t$, assumption \eqref{eq:H07} guarantees that  
$\Gamma \subset \partial V(t)\cap \partial \Omega$ for all
$t>0$.  In fact,  for any $x \in \Gamma$ one finds  $x_j \in
V_0\subset V(t)$ such that $x_j \to x$. As $x \not \in V(t)$, this
implies that   $ x
\in  \partial V(t)\cap \partial \Omega$. Moreover, one can find
$\gamma_0\in (0,\ove \epsi)$ such that
$ \inf_{g\in \Gamma}  d(g,\partial V (t) \setminus \partial
\Omega)\geq \gamma_0$  for all $t>0$. 

As already mentioned, the open set $V(t)$ may not be smooth, see Proposition~\ref{prop:john}. Nonetheless, under
assumptions \eqref{eq:H05} and \eqref{eq:H07}, the
portion $\Gamma$ of $\partial V(t)$ is smooth and is well
separated from $\partial V(t) \setminus \partial \Omega$ for all $t>0$. This allows us to
define the trace  on $\Gamma$  of any $u \in H^1(V(t))$. More
precisely, set  $G\subset (\Gamma + B_\delta) \cap \Omega$ 
for some $ 0 <  \delta< \inf_{g\in \Gamma}
d(g,\partial V_0\setminus \partial \Omega) $
small enough with $G  \subset V_0  $ Lipschitz  and
$\Gamma \subset \partial G$.  Then, the trace operator on $H^1(V(t))$ can
be defined by composing the restriction $u \in H^1(V(t))
\mapsto u|_G\in  H^1(G)$ and the
 trace operator $\gamma:H^1(G) \to H^{1/2}(\Gamma)$. This allows us to
 specify 
$$H^1_\Gamma(V(t)):=\{u \in H^1(V(t)) \ | \ u =0 \ \text{on}
\ \Gamma\}\quad \text{for} \ t >0.$$
 Note that the trace operator $\gamma$ depends on $\Gamma$ but not
on the choice of $G$. 
 
The Hamiltonian $H$ is asked to fulfill the following
 \begin{align}
    &H : \ove \Omega \times \Rz \times \Rz^n\to\Rz \quad \text{is continuous},\label{eq:nH01u}\\
    &C_{xu}:= \{p\in \Rz^n \ | \ H(x,u,p)\leq 0\} \quad  \text{is convex}
      \quad 
      \forall (x,u)\in \ove \Omega \times \Rz, \label{eq:nH02u}\\ 
     &\forall (x,u)\in \ove \Omega \times \Rz \ \forall p\in C_{xu} \ \forall \delta >0 \
       \exists p_\delta \in \Rz^n: \ H(x,u,p_\delta)<0 \ \text{and} \
       |p-p_\delta|<\delta, \label{eq:H03u}\\
   & \forall (x,u,p)\in \ove \Omega \times \Rz \times \Rz^n: \
      H(x,u,p)< 0 \ \Leftrightarrow \exists
     \lambda \in (0,1) \ \text{such that} \  p/\lambda \in C_{xu},  \label{eq:H03ub}\\
  & \exists 0<\sigma_* \leq \sigma^*: \quad B_{\sigma_*}\subset
   C_{xu}\subset   \ove{B_{\sigma^*}}\quad \forall (x,u)\in \ove \Omega
    \times \Rz.\label{eq:nH04u} 
 \end{align}
 Note that $p\mapsto H(x,u,p)$ is not asked to be convex but merely to
 have convex sublevels with respect to $p$.  Assumption
 \eqref{eq:nH04u}  ensures that $0 \in C_{xu}$ for all  $ (x,u)\in \ove
 \Omega \times \Rz$, so that \eqref{eq:H03ub} gives $H(x,u,0)<0$.  The
 simplest setting where \eqref{eq:nH01u}--\eqref{eq:nH04u} hold is that
 of the {\it generalized eikonal equation}, namely, 
 $H(x,u,p) = \gamma(x,u)|p|-1$ with $\gamma$  uniformly  continuous and
 $1/\sigma^* \leq \gamma(\cdot ) \leq 1/\sigma_*$

For all $(x,u)\in \ove \Omega \times \Rz$ we define the  {\it support function}
$q\in \Rz^n\mapsto\sigma(x,u,q)$ of the nonempty, convex, and closed set $C_{xu}$ as 
\begin{equation*}
  \sigma(x,u,q) = \sup \{q\cdot p \ | \ H(x,u,p)\leq
0\}. 
\end{equation*}
 The functions  $q \in \Rz^n\mapsto \sigma(x,u,q)$ are   nonnegative and positively
$1$-homogeneous.  In particular,  $0 = \sigma(x,u,0)=\min_q \sigma(x,u,q)$.

 The relaxation kernels are asked to fulfill
 \begin{equation}
  \label{eq:H06}
  k \in W^{1,1}(\Rz_+), \quad\phi\in
H^1(\Rz^n).
\end{equation}
Given $v
\in C(\ove \Omega)$, we set $V(t)=\{x \in  \Omega \ | \ v(x) <t\}$ and
$Q_v:=\cup_{t>0}V(t)\times \{t\}$. 
Let us note that $Q_v=\{(x,t)\in  \Omega
\times \Rz_+ \ | \ v(x)-t< 0\}$ with $(x,t)\in \ove \Omega
\times \Rz_+\mapsto v(x)-t$ continuous. In particular, $Q_v\subset
 \Rz^n  \times \Rz_+$ is measurable.

 For all $z \in L^\infty_{\rm loc}( [0,\infty);L^2(\Rz^n))$ we  define the space-time convolution  
\begin{equation*} K
   z  (x,t):=\int_0^t\!\!\int_{\Rz^{ n}}k(t-s)\phi(x-y)\, z  (y,s)\,
\d y\,  \d s.
\end{equation*}
 Young's convolution inequality 
ensures  that $Kz$ is Lipschitz continuous in space-time. More
precisely,  from
\begin{align*}
&  \nabla (Kz)(x,t) = 
\int_0^t\!\!\int_{\Rz^n}k(t-s)\nabla \phi(x-y) \,z (y,s)\, \d y \,
  \d s,\\
  & (Kz)_t(x,t) = k(0)\int_{\Rz^n}\phi(x-y)\,z(y,t) \,
\d y
+\int_0^t\!\!\int_{\Rz^n}k'(t-s) \phi(x-y) \,z (y,s)\, \d y \,
\d s,
\end{align*}
we readily check that 
$$\forall T>0 \ \exists c(T)>0: \  \| Kz \|_{ W^{1,\infty}(\Rz^n \times  (0,T)) }\leq c(T)
  \|z\|_{L^\infty(0,T;L^2(\Rz^n))}$$
where the constant $c(T)$ depends also  on $\| k \|_{ W^{1,1}(\Rz_+)}$ and
$\|  \phi\|_{
  H^1(\Rz^{ n})}$. 
 For $u:Q_v \to \Rz$ 
measurable with $ t \mapsto \| u(\cdot, t) \|_{L^2(V(t))}\in
L^\infty_{\rm loc}( [0,\infty) )$ we let $Ku:=K\ove u$, where $\ove u(\cdot,t) $ denotes the trivial extension of $u(\cdot,t)$
to the whole $\Rz^n$ for a.e. $t>0$.

We are now ready to state our main result. 

\begin{theorem}[Existence for the free boundary problem]\label{thm:main}
  Under assumptions \eqref{eq:H05}--\eqref{eq:H06} there  exist
  $v \in C(\ove \Omega)$ and   $u:Q_v \to \Rz$
  measurable  with  $u(\cdot,t) \in H^1_{\Gamma}(V(t))$ for a.e. $t
  >0$ such that  $t  \in [0,\infty)  \mapsto \| u(\cdot, t)\|_{H^1(V(t))}$ is
  locally bounded,  
\begin{align}
    &\int_{V(t)}\nabla u \cdot \nabla w \, \d x =
    \int_{V(t)}  w \, \d x\quad \forall w \in
      H^1_\Gamma(\Omega), \ \text{for a.e.} \ t >0, \label{eq:0111}\\
&v(x) = \min \Bigg\{\int_0^{1} 
  \sigma\big(\gamma(s),Ku\big(\gamma(s),v(\gamma(s))\big),\gamma'(s)\big) \ {\d s} \ \Big| \ \gamma \in  {\mathcal C}_x , \ \gamma (\cdot)\in
                                    \ove \Omega\Bigg\}\quad \forall x
                  \in \ove \Omega,  \label{eq:0222}
 \end{align}
 and $v$ is a \emph{viscosity solution} to
\eqref{eq:02}, namely, for any $\varphi \in C^1(\Rz^n)$, we have that $$H(x,Ku(x,v(x)),\nabla
\varphi(x))\leq (\geq) \, 0$$ at any local minimum (maximum, respectively) point $x\in \Omega\setminus \ove
{V_0}$ of $\varphi - v$.

Moreover,  the function $v$ is a
\emph{constrained viscosity solution}, namely, for any $\varphi \in C^1(\Rz^n)$, we have that $H(x,Ku(x,v(x)),\nabla
\varphi(x))\geq 0$ at any local  maximum  point $x\in \partial \Omega\setminus \ove
{V_0}$ of $\varphi - v$.

\end{theorem}

We remark  in particular  that the minimum in
\eqref{eq:0222} is always attained. Note that the function $v$ appears in the argument
of the integrand in \eqref{eq:0222}, as well. As such, \eqref{eq:0222} does not directly
define $v$, given $u$. In fact, 
\eqref{eq:0222}  does not imply uniqueness,  in general.  Moreover,  note that the test function $w $
in \eqref{eq:0111} is taken in the whole $H^1_\Gamma(\Omega)$ instead of
$H^1_\Gamma(V(t))$ at $t>0$. Given $w \in H^1_\Gamma(\Omega)$, its restriction
$w|_{V(t)}$ belongs to $H^1_\Gamma(V(t))$  for all $t
>0$. On the other hand, $V(t)$ need not be an extension domain,
 see Figure~\ref{fig:1}, 
so that a function in  $ H^1_\Gamma(V(t))$ may have no extension to
$H^1_\Gamma(\Omega)$. Hence, satisfying  \eqref{eq:0111} for $w \in H^1_\Gamma(\Omega)$
delivers a weaker statement  than requiring    it   for  $w\in
H^1_\Gamma(V(t))$ at $t>0$. However,
\eqref{eq:0111} still implies that  $u$ solves the elliptic problem in
 the 
distributional sense, and  that  the homogeneous Neumann condition at
$\partial V(t)\setminus \ove \Gamma  $  is locally  recovered,   wherever the
boundary $ \partial V(t)$ is smooth enough  and $V(t)$ lies
locally on one side of its boundary.

Theorem~\ref{thm:main} is proved in Section \ref{sec:elastic} via a
 fixed-point argument.

\section{Some remarks on the Hamilton--Jacobi subproblem}\label{sec:setting}

In preparation for the proof of Theorem~\ref{thm:main}, we collect some material on the Hamilton--Jacobi
subproblem. More precisely, given $\Omega \subset \Rz^n$ fulfilling
\eqref{eq:H05} and 
 \begin{equation}\label{eq:H00}
   V_0 \subset \Omega \quad \text{nonempty and open},
     \end{equation}
in this section, we are interested in the solvability in the viscosity
sense of problem
 \begin{align}
  &H(x,\nabla v )=0\quad \text{in} \ \Omega\setminus
    \ove{V_0}, \label{eq:1}\\
   & v =0 \quad \text{in} \ V_0.\label{eq:2}
 \end{align}
 
  We recall that $v\in C(\ove \Omega)$ is said to be
  a {\it viscosity  solution} to equation
  \eqref{eq:1} if, for any $\varphi \in C^1(\Rz^n)$, we have that $H(x,\nabla
\varphi(x))\leq (\geq)\, 0$ at any local minimum (maximum, respectively) point $x\in \Omega\setminus \ove
{V_0}$ of $\varphi - v$.  
In case $\Omega=\Rz^n$, we refer to problem \eqref{eq:1}--\eqref{eq:2}  
as {\it unconstrained}, whereas, if $\Omega \not = \Rz^n$ we say that
problem \eqref{eq:1}--\eqref{eq:2} is {\it constrained}. 
A viscosity solution $v\in C(\ove \Omega)$ to \eqref{eq:1} is said to be
  a {\it constrained  viscosity  solution} if 
  \begin{equation*} 
      H(x,\nabla v)\geq 0\quad \text{on} \ \partial \Omega\setminus
    \ove{V_0}
  \end{equation*}
  in the viscosity sense, namely, for any $\varphi \in C^1(\Rz^n)$, we have that $H(x,\nabla
\varphi(x))\geq 0$ at any local  maximum  point $x\in \partial \Omega\setminus \ove
{V_0}$ of $\varphi - v$.

We assume that 
\begin{align}
  & H: \ove \Omega\times \Rz^n\to\Rz \  \text{is continuous,} \  
    C:=\{(x,p) \in \ove \Omega \times \Rz^n\ | \ H(x,p)\leq 0\}, \label{eq:H01}\\
    &C_x:= \{p\in \Rz^n \ | \ H(x,p)\leq 0\} \quad  \text{is convex}
      \quad 
      \forall x \in \ove \Omega, \label{eq:H02}\\
    &\forall (x,p) \in C\ \forall \delta>0 \ \exists (x,p_\delta) \in
     {\rm int}_{\ove \Omega \times \Rz^n} {C} \ \text{such that} \ |p-p_\delta|\leq
      \delta,  \label{eq:H03}\\
    & \forall (x,p)\in  \ove \Omega \times \Rz^n: \ H(x,p)<0 \ \Leftrightarrow \ \exists
      \lambda\in (0,1):\  \   p/\lambda \in C_x,  \label{eq:H03b}\\
  & \exists 0<\sigma_* \leq \sigma^*: \quad B_{\sigma_*}\subset
   C_x\subset   \ove{B_{\sigma^*}}\quad \forall x \in \ove \Omega. \label{eq:H04}
  \end{align}
  Note that these assumptions  correspond to 
  \eqref{eq:nH01u}--\eqref{eq:nH04u}.  Again, assumptions \eqref{eq:H03b}--\eqref{eq:H04}
  ensure that $H(x,0)<0$ for all $x \in \ove \Omega$. 

Due to the
definition of viscosity solution (or constrained viscosity solution), in order to solve problem   \eqref{eq:1}--\eqref{eq:2} we can equivalently replace  $H$
by any $\ove H$ such that $\{\ove H\leq  0\}= \{ H\leq 0\}$ and
$\{\ove H\geq  0\}= \{ H\geq 0\}$. This fact is used in the proof
of Proposition~\ref{prop:representation}  below.

For all $x\in \ove \Omega$ we denote by
  $\sigma(x,\cdot):\Rz^n \to [0,\infty)$ the support function of the nonempty,
convex, and closed set $C_x$, namely,
$$\sigma(x,q) = \sup\{q\cdot p\ | \ H(x,p) \leq 0\}\quad \forall
(x,q)\in \ove \Omega \times \Rz^n. $$
We have that $q\mapsto \sigma(x,q) $ is nonnegative and positively
$1$-homogeneous.  In particular,   $0 = \sigma(x,0)=\min_q \sigma(x,q)$. Assumption \eqref{eq:H04} ensures that 
\begin{equation}
  \sigma_* |q| \leq \sigma(x,q) \leq  \sigma^*|q|\quad \forall
(x,q)\in \ove \Omega \times \Rz^n.\label{eq:sigma}
\end{equation}
Let us collect some basic properties of $\sigma$ in the following.

\begin{lemma}[Properties of $\sigma$]\label{lemma:sigma}
 Under assumptions \eqref{eq:H01}--\eqref{eq:H03}  and   \eqref{eq:H04}, the
 support function $\sigma$ is
 continuous and $q\mapsto \sigma(x,q)$ is convex for all $x\in \ove \Omega$.
\end{lemma}
\begin{proof}
  For all $x\in \ove \Omega$, the function $q\mapsto \sigma(x,q)$ is the
  support function of the convex set $C_x$  and   is convex. 
  
  Let $(x_j,q_j)_j \subset  \ove \Omega \times \Rz^n$ with $(x_j,q_j) \to (x,q)$ and, for all $\zeta>0$ let $p_\zeta \in C_x$
such that $\sigma(x,q)\leq q\cdot p_\zeta
+\zeta$. Using \eqref{eq:H03}  we find
$ p_{\zeta\delta} $ such that $(x,p_{\zeta\delta}) \in C^\circ$ and $q\cdot p_\zeta \leq
q\cdot  p_{\zeta\delta} +|q|\delta$. In particular, we have that 
$(x_j, p_{\zeta\delta})\in C$ for $j$ large enough,
so that $ p_{\zeta\delta} \in
C_{x_j}$, and we conclude that
\begin{align*}
  &\sigma(x,q)\leq q\cdot p_{\zeta}+\zeta \leq  q\cdot  p_{\zeta\delta} +
\zeta+|q|\delta  
    = \lim_{j\to \infty}q_j\cdot p_{\zeta\delta}  + \zeta+|q|\delta \\
  &\quad \stackrel{ p_{\zeta\delta}\in
C_{x_j}}{\leq }\liminf_{j\to
    \infty}\sigma(x_j,q_j) +\zeta+|q|\delta  .
    \end{align*}
As $\zeta$ and $\delta $ are arbitrary, this proves that $\sigma$ is lower
semicontinuous.

Let us now check upper semicontinuity. For all $j$ we
can find $p_j\in C_{x_j}$ such that $ 
\sigma(x_j,q_j)\leq q_j \cdot p_j + 1/j$. Let the subsequence
$(j_k)_k$ be such that $
\lim_{k\to \infty}\sigma(x_{j_k},q_{j_k})=\limsup_{j\to\infty}\sigma(x_j,q_j)$. As $p_{j_k}\in
C_{x_{j_k}}\subset B_{\sigma^*}$ by \eqref{eq:H04}, one can extract
without relabeling in such a way that $p_{j_k}\to p$. As $C$ is closed
by \eqref{eq:H01}, we have that $(x,p) \in C$, namely, $p\in
C_x$, and  
$$ \limsup_{ j\to\infty}\sigma(x_j,q_j)=
\lim_{k\to \infty}\sigma(x_{j_k},q_{j_k}) = \lim_{k\to \infty}q_{j_k}
\cdot p_{j_k}  = q\cdot p  \stackrel{p\in C_x}{\leq}
\sigma(x,q) $$
so that $\sigma$ is upper semicontinuous, as well. 
\end{proof}

In the unconstrained case $\Omega=\Rz^n$, existence of viscosity solutions
to \eqref{eq:1}--\eqref{eq:2} is ensured by the following. 
  
  \begin{proposition}[Existence and representation in the
    unconstrained case]\label{prop:rn}
   Let $\Omega=\Rz^n$,  assume \eqref{eq:H02}--\eqref{eq:H04},
   and let $C=\{(x,p)\in \Rz^n \times \Rz^n \ | \
   H(x,p)\leq 0\}$ be closed. Then,   problem
    \eqref{eq:1}--\eqref{eq:2} admits a unique nonnegative viscosity
    solution, which is given by the representation formula
    \begin{align}
      v(x) = \min \left\{\int_0^{1} 
      \sigma(\gamma(s),\gamma'(s)) \ {\d s} \ \Big| \  \gamma \in  {\mathcal C}_x \right\}\quad  \forall x
      \in \Rz^n . \label{eq:vis3}
    \end{align}
  \end{proposition}
  \begin{proof} We first check that 
    \begin{equation}
      \label{eq:Mennucci}
C_x = \overline{\{p\in \Rz^n \ | \ (x,p)\in C^\circ\}}=: D_x \quad \forall x \in \Rz^n.
\end{equation}
The inclusion $D_x\subset C_x$ is immediate: for
all 
$p\in D_x$ we find $p_j \to p$ with $H(x,p_j)\leq 0$. As $C\ni (x,p_j)
\to (x,p)$, the  closedness  of $C$ ensures that $(x,p)\in C$, namely,  $p
\in C_x$.

Let now $p \in C_x$. Condition \eqref{eq:H03} implies that we can find
a sequence $(p_\delta)_\delta $ with $(x,p_\delta)\in C^\circ$ and
$p_\delta \to p$ as $\delta \to 0$. This proves that $p\in D_x$, hence $C_x
\subset D_x$.

Having proved \eqref{eq:Mennucci} and  recalling that $H(x,0)<0$ for all
$x \in \Rz^n$  the  assertion  follows from
    \cite[Thm.~3.15]{Mennucci}.
  \end{proof}
  
A consequence of Proposition~\ref{prop:rn} is that the  
minimum in \eqref{eq:vis3} is always attained. In 
 the following, we call {\it optimal curve} for $x$ a minimizer
 $\gamma_x\in  {\mathcal C}_x $ in 
\eqref{eq:vis3}. Note that optimal curves $\gamma_x$
are not necessarily unique.  Let   us mention that the assumptions of
Proposition~\ref{prop:rn} could be relaxed (see \cite{Mennucci}). We
work nonetheless in the setting of  \eqref{eq:H01}--\eqref{eq:H04}, which  is tailored to   the coupling with the PDE for the
activating field $u$.

By dropping the nonnegativity requirement on
$v$, uniqueness in Proposition~\ref{prop:rn} fails. An
elementary illustration of this fact  is given by   the viscosity solutions 
\begin{equation}
    x\mapsto v_0(x) =
(|x|-1)^+\quad \text{and} \quad x\mapsto v_r(x) =
\min\{v_0(x), r-|x|\}\quad \forall r \geq 1,\label{eq:count}
\end{equation}
to problem \eqref{eq:1}--\eqref{eq:2} in the
 unconstrained case $\Omega=\Rz^n$ for
 $H(x,p)=|p|-1$ and  
  $V_0 ={B_1}$.  
 We now turn to the constrained case $\Omega \not = \Rz^n$
 under assumption \eqref{eq:H05}. 
In contrast  to   the unconstrained case, uniqueness of viscosity
solutions to the constrained case cannot be restored by simply requiring
nonnegativity: Consider again the case of 
 \eqref{eq:count} and take $\Omega=B_3$. The restrictions to $\Omega$ of all 
 functions in \eqref{eq:count} are viscosity solutions to the
 constrained problem.  Moreover,   we have that $v_0$, $v_r\geq 0$ in $\Omega$
 for  all  $r  \geq 3$, so that nonnegativity alone is not
 selective.   The   unique constrained viscosity solution in $\Omega =
 B_3$ is $v_0$, coinciding indeed with $v_r$ for all $r \geq 5$.

Aiming at identifying some possible selection
principle among viscosity solutions, we are interested in a
constrained version of the representation formula \eqref{eq:vis3}.
The existence result in the constrained case reads as follows.

\begin{proposition}[Existence and representation in the constrained
  case]\label{prop:representation} Assume \eqref{eq:H05}--\eqref{eq:E2},
  \eqref{eq:H00}, and
  \eqref{eq:H01}--\eqref{eq:H04}. Then, $v \in
  C(\ove \Omega)$ given by
\begin{align} 
     v(x) = \min \Bigg\{\int_0^{1} 
  \sigma(\gamma(s),\gamma'(s)) \ {\d s} \  \Big|  \  &\gamma \in  {\mathcal C}_x , \ \gamma (\cdot)\in
                                    \ove \Omega\Bigg\}\quad \forall x
                  \in \ove \Omega.  \label{eq:vis4}
\end{align} 
is a  constrained  viscosity
 solution to  
 problem \eqref{eq:1}--\eqref{eq:2}. 
  \end{proposition}
  
  \begin{remark}[Alternative assumptions]\label{rem:prima}
     If $p\in \Rz^n \mapsto H(x,p)$ is convex and $\inf
    H(x,\cdot)<0$ for all $ x\in \ove \Omega$, property \eqref{eq:H03}
    follows. Moreover, if $\Omega $ is convex, assumption \eqref{eq:H03b} can be
    dropped. This is also the case if $\ove \Omega$ is
  homeomorphic to a convex and closed set $D\subset \Rz^n$ via a
  global homeomorphism
  $h: \Rz^n \to \Rz^n$, namely, $\ove \Omega=h(D) $, see
  Remark~\ref{rem:seconda} below.
\end{remark}

\begin{remark}[Uniqueness]\label{rem:uniqueness} If $p\mapsto H(x,p)$ is not convex,   uniqueness
  cannot 
  be expected in general, see \cite{Lions}. On the other hand, if 
$H$
  is Lipschitz continuous  and 
  $p\mapsto H(x,p)$ is convex  with  $H(x,0)<0$ for all $x \in \ove\Omega$, \cite[Prop.~IX.3]{CapuzzoDolcetta-Lions} guarantees that
  $v$ given by \eqref{eq:vis4} is the unique  constrained 
  viscosity solution to  problem \eqref{eq:1}--\eqref{eq:2}. This in
  particular applies to the case of the generalized eikonal equation
  $H(x,p)=\gamma(x)|p|-1$ for $\gamma$ Lipschitz continuous, bounded,
  and uniformly positive.
\end{remark}

\begin{proof}[Proof of Proposition~\ref{prop:representation}] The proof relies on a penalization procedure  in the
  spirit of \cite[Sec.~VII]{CapuzzoDolcetta-Lions}, and is
  divided into steps. At first, we reformulate the problem by
  equivalently replacing $H$ by another Hamiltonian $\ove H$ having  $\{\ove H \leq 0\}=
\{H\leq 0\}$ and  $\{\ove H \geq 0\}=
\{H\geq 0\}$  (Step 1). Then, we
  provide extensions $\ove H_\epsi$ of $\ove H$ to all $x \in \Rz^n$, depending on a small parameter
  $\epsi>0$ (Step 2). For all such $\ove H_\epsi$ we solve the unconstrained
  problem in $\Rz^n$ by the representation formula \eqref{eq:vis3}. This defines
  functions $v_\epsi\in C(\Rz^n)$ (Step 3). By taking $\epsi \to 0$, we
  prove that $v_\epsi \to v$ locally uniformly in $\ove \Omega$, that
  $v\in C(\ove \Omega)$ is a viscosity solution to
  \eqref{eq:1}--\eqref{eq:2} (Step 4), and that the representation formula
  \eqref{eq:vis4} holds (Step 5).  Finally,    we prove that $v$ is a constrained viscosity
  solution (Step 6). 
  
  \medskip
  \noindent 
  {\it Step 1: Replacing $H$ by $\ove H$.} Define
  \begin{equation*}
  \ove H(x,p):= \inf\{\lambda >0 \ | \ p/\lambda \in
C_x\}-1\quad \forall (x,p)\in \ove \Omega \times \Rz^n.
\end{equation*}
For all $x \in \ove \Omega$, $\ove H(x,\cdot)$ is the standard {\it Minkowski functional} of the convex set
$C_x$, up to an
additive constant. As $0\in B_{\sigma_*}\subset C_x$ we have that 
$\ove H(x,0)=-1$ and $p\mapsto \ove H(x,p)+1$  is positively
$1$-homogeneous.  In particular, owing to \eqref{eq:H04} we have
that
$$\ove H(x,p)+1 \leq \frac{|p|}{\sigma_*} \quad \forall (x,p)\in \ove
\Omega\times \Rz^n$$
so that $\ove H$ is finite. 
Moreover, for all $(x,p)\in \ove \Omega \times \Rz^n$ with $H(x,p)\leq 0$ one has
that $(x,p/1)\in C$, hence $\ove H(x,p)\leq 0$. On the other hand, if
$\ove H(x,p)\leq 0$ there exists $\lambda_\delta \to  \lambda \leq
 1$ such that
$(x,p/\lambda_\delta)\in C$ and the  closedness  of $C$ from \eqref{eq:H01}
ensures that $(x, p/\lambda)\in
C$.  As $C_x$ is convex and $0 \in C_x$ we have that $(x,p) \in C$,  namely, $H(x,p)\leq 0$. This proves that $\{\ove H \leq 0\}=
\{H\leq 0\}$. As assumption \eqref{eq:H03b} directly implies
that $\{H<0\}=\{\ove H<0\}$, by taking the complements we have
that $ \{\ove H \geq 0\}=\{H
    \geq 0\}$. We are hence allowed to replace  $H$ by $\ove H$ in
  \eqref{eq:1}--\eqref{eq:2} without changing the corresponding
  viscosity solutions.  
  Note more precisely that 
 $C_x=\{p \in\Rz^n
  \ | \ \ove H(x,p)\leq
  0\}$ for all $x \in \ove \Omega$, so that $\sigma$ and formula
  \eqref{eq:vis4} remain unchanged.

   We now prove that
  \begin{equation}\label{eq:Hcontinuous}
    \ove H \ \ \text{is continuous.}
  \end{equation}
  At first, using
  \eqref{eq:H03b} we readily check that for all $t>0$ and $(x,p)\in
  \ove \Omega \times \Rz^n$
  \begin{align}
    &\ove H(x,p)+1<t \ \Leftrightarrow \ H(x,p/t)<0, \label{eq:sublevels1}\\
    &\ove H(x,p)+1\leq t \ \Leftrightarrow \ H(x,p/t)\leq 0. \label{eq:sublevels2}
  \end{align}
  Let us check \eqref{eq:sublevels1}: As $\ove H(x,\cdot) +1 $ is positively
  $1$-homogeneous, one has that $ \ove H(x,p)+1<t $ is equivalent to
  $\ove H(x,p/t)+1<1$, which in turn is equivalent to the existence of
  $\lambda \in (0,1)$ such that $p/(\lambda t)\in C_x$, which by
  \eqref{eq:H03b} 
  is equivalent to $H(x,p/t)<0$. Equivalence \eqref{eq:sublevels2} is
  proved analogously.

  Let now $(x_j,p_j)\to (x,p)$. Take $t>\ove H(x,p)+1$ and use
  \eqref{eq:sublevels1} to get that $H(x,p/t)<0$. As $H$ is
  continuous,  $H(x_j,p_j/t)<0$ for $j$ large enough. Again
  \eqref{eq:sublevels1} gives that $\ove H(x_j,p_j)+1<t$ for $j$
large enough. Hence, $\limsup_{j\to \infty}\ove H(x_j,p_j)+1\leq t$ for all
$t>\ove H(x,p)+1$. This implies that  $\limsup_{j\to \infty}\ove
H(x_j,p_j)\leq \ove H(x,p)$, so that $\ove H$ is upper semicontinuous.

Extract now without relabeling in such a way that $\ove H(x_j,p_j) \to
\liminf_{k\to \infty}\ove H(x_k,p_k)=:\iota$. For all $t >\iota+1$ we
have that $\ove H(x_j,p_j)+1<t$ for $j$
large enough. Using \eqref{eq:sublevels1} this ensures that $H(x_j,p_j/t)<0$ for $j$
large enough. By the continuity of $H$ we hence have that $H(x,p/t)\leq
0$ and \eqref{eq:sublevels2} ensures that $\ove H(x,p)+1\leq t$. As $t
>\iota+1$ is arbitrary, we have proved that $\ove H(x,p) \leq \iota =
\liminf_{k\to \infty}\ove H(x_k,p_k)$, so that $\ove H$ is lower semicontinuous.

  \medskip
  \noindent 
  {\it Step 2: Extension of $\ove H$.}  For all $\epsi \in (0,\min\{\sqrt{\ove \epsi},1/\sigma^*\})$
  we define
  $$C_{x\epsi}:=\left(1-\frac{d(x,\Omega)}{\epsi^2} \right)^+C_{
    \ove \pi  (x)} +
  \min\left\{ \frac{d(x,\Omega)}{\epsi^2},1\right\}\ove{B_{1/\epsi}}\quad
  \forall x \in \Rz^n.$$
  In the above right-hand side, the first set  reduces to $\{0\}$ 
  unless $d(x,\Omega)<\epsi^2$. In this case, as $\epsi^2 <\ove \epsi$ we have that
  $d(x,\Omega)<\ove \epsi$ and the projection $ \ove \pi : \ove
  \Omega  
  + B_{\ove \epsi} \to \ove \Omega $ given by $\ove \pi  (x) = {\rm
    arg\,min}\{|x- \omega  |\:|\:\omega  \in \ove
\Omega\}$ is
  well defined. Specifically, we have $C_{x\epsi} = C_x$ for all $x
  \in \ove \Omega$ and $C_{x\epsi}=\ove{B_{1/\epsi}}$ if $d(x,\Omega)\geq
  \epsi^2$.

  Let us now define $\ove H_\epsi :\Rz^n \times \Rz^n \to [ -1 ,\infty)$
  as
  $$\ove H_\epsi(x,p) = \inf\{\lambda>0 \ | \ p/\lambda \in
  C_{x\epsi}\}-1$$
  and observe that $C_{x\epsi}=\{p\in \Rz^n \:|\: \ove
  H_\epsi(x,p)\leq 0\}$ for all $x \in \Rz^n$.

  In the remainder of this step, we  check that assumptions  \eqref{eq:H01}--\eqref{eq:H04} for $H$ in $\ove \Omega$ imply
 that $ \ove H_\epsi $ fulfills  the 
  assumptions of Proposition \ref{prop:rn}. 
  To simplify notation, in the following we use  $\alpha(x):=
  (1-d(x,\Omega)/\epsi^2)^+ $ and $\beta(x):=\min\{d(x,\Omega)/\epsi^2,1\}$ for
  all $x \in \Rz^n$.

  At first, let us check that $C_\epsi:=\{(x,p)\in \Rz^n \times \Rz^n
  \ | \ \ove H_\epsi(x,p)\leq 0\}$ is closed. Take a sequence
  $(x_j,p_j)\in C_\epsi$ with $(x_j,p_j)\to (x,p)$. Then, we can find
    $q_j\in C_{ \ove\pi (x_j)}$ and $r_j\in \ove{B_{1/\epsi}}$ such that
    $p_j = \alpha(x_j)q_j + \beta(x_j)r_j$.
 Upon possibly extracting not relabeled subsequences and using
    \eqref{eq:H01} and \eqref{eq:H04} we can assume that $q_j \to q \in C_{ \ove\pi (x)}$ and
    $r_j \to r \in \ove{B_{1/\epsi}}$. Hence,
    $p_j \to \alpha(x)q + \beta(x) r \in
    C_{x\epsi}$. This proves that $(x,p)\in C_\epsi$. In particular,
    $C_\epsi$ is closed.
    
   For all $x\in\Rz^n$, $C_{x\epsi}$ is the sum of two convex sets and
  hence convex. This implies that \eqref{eq:H02} holds for
  $\ove H_\epsi$ for all $x \in \Rz^n$. 

  In order to check \eqref{eq:H03} for $\ove H_\epsi$, we fix $(x,p)\in
  C_\epsi$ and distinguish some cases:
  \begin{itemize}
  \item Case $x \in \Omega$. From $p\in  C_{x\epsi}= C_x$, for all $\delta
    >0$ we use \eqref{eq:H03}
    for $H$ in order to find $p_\delta$ with $|p-p_\delta|\leq \delta$
    and $(x,p_\delta)\in  {\rm int}_{\ove \Omega \times \Rz^n}C$. As $\Omega$ is open, we directly get
    that $(x,p_\delta)\in C_\epsi^\circ$, as well.
    \item Case  $x \in \partial \Omega$.  For all $p \in
      C_{x\epsi}=C_x$ and all $\delta>0$ we again use \eqref{eq:H03}
      and find $p_\delta$ such that $|p-p_\delta|\leq \delta$ and
      $(x,p_\delta)\in  {\rm int}_{\ove \Omega \times \Rz^n}C$. Let $(\widetilde x,\widetilde p)$ be close to
      $(x,p_\delta)$, in particular let $y =\widetilde x - x$ and $s=\widetilde p - p_\delta$ be
      small. Then
      \begin{align*}
    \alpha^{-1}(\widetilde x) \widetilde p &= p_\delta + s + (\alpha^{-1}(\widetilde
                                    x){-}1)\widetilde
    p=p_\delta + s + (\alpha^{-1}(\widetilde x){-}1)(p_\delta +s)
      \end{align*}
       and $\alpha^{-1}(\widetilde x) \widetilde p $ is close to
      $p_\delta $  for $y$ and $s$ small, as $\lim_{y\to
        0}(\alpha^{-1}(\widetilde x)-1)=0$. Hence, $( \ove\pi  (\widetilde x),
      \alpha^{-1}(\widetilde x) \widetilde p)$ is close to
      $(x,p_\delta) \in  {\rm int}_{\ove \Omega \times \Rz^n}C$, 
      so that $ \alpha^{-1}(\widetilde x) \widetilde p \in
      C_{ \ove\pi (\widetilde x)}$.  This proves that
      $\widetilde p = \alpha(\widetilde x) (  \alpha^{-1}(\widetilde x) \widetilde p
      )\in C_{\widetilde x\epsi}$, hence $(\widetilde x, \widetilde p )\in
      C_\epsi$. We have checked that $(x,p_\delta)\in C_\epsi^\circ$.
  \item Case $0<d(x,\Omega)<\epsi^2$. We have $p=\alpha(x)q+\beta(x)r$ with
    $q \in C_{ \ove\pi (x)}$, $|r|\leq 1/\epsi$, and
    $\alpha(x),\,\beta(x)>0$, hence $x\mapsto 1/\alpha(x)$ and
    $x\mapsto 1/\beta(x)$ are  well-defined and are continuous in a neighborhood
    of $x$. For all $\delta>0$, using \eqref{eq:H03} we find
    $q_\delta$ with $|q-q_\delta|\leq \alpha^{-1}(x)\delta/2$ and
    $( \ove\pi (x),q_\delta)\in  {\rm int}_{\ove \Omega \times \Rz^n}C$. Moreover, we find $r_\delta$ with
    $|r-r_\delta|\leq \beta^{-1}(x) \delta/2$ and
    $|r_\delta|<1/\epsi$. Set now $p_\delta:=\alpha(x)q_\delta +
    \beta(x) r_\delta$, so that $|p-p_\delta|\leq \delta$. Let $(\widetilde x,\widetilde p)$ be close to
    $(x,p_\delta)$. Setting $y =\widetilde x - x$ and $s=\widetilde p - p_\delta$ one has that
    \begin{align*}
      \widetilde p &= \alpha(\widetilde x)\left(q_\delta {+}
      \left(\frac{\alpha(x)}{\alpha(\widetilde x)}{-}1\right)q_\delta
        {+}\frac{s}{2\alpha(\widetilde x)}\right)  + \beta(\widetilde x)\left(r_\delta {+}
      \left(\frac{\beta(x)}{\beta(\widetilde x)}{-}1\right)r_\delta
        {+}\frac{s}{2\beta(\widetilde x)} \right)\\
      & =: \alpha(\widetilde x) \haz q +\beta(\widetilde x)\haz r .
    \end{align*}
For $y$ and $s$ small, $\haz q$ and $\haz r$ are arbitrarily close to $q_\delta$
and $r_\delta$, respectively. In particular,  we have that $\haz q
\in C_{\ove \pi(\widetilde x)}$ and  $\haz r\in
\ove{B_{1/\epsi}}$. This proves that $\widetilde p\in C_{\widetilde
  x\epsi}$ or, equivalently, $(\widetilde x,\widetilde p) \in C_\epsi$. As a
consequence, we have that $(x,p_\delta)\in C_\epsi^\circ$.
\item Case $d(x,\Omega)=\epsi^2$. As $p\in \ove{B_{1/\epsi}}$, for all
  $\delta >0$ we can
  take $p_\delta \in B_{1/\epsi}$ with $|p-p_\delta|\leq
  \delta$. Let $(\widetilde x,\widetilde p) $ be close to $(x,p_\delta)$. In
  particular, let $y =\widetilde x - x$ and $s=\widetilde p - p_\delta$ be small. Then
  \begin{align*}
    \beta^{-1}(\widetilde x) \widetilde p &= p_\delta + s + (\beta^{-1}(\widetilde
                                    x){-}1)\widetilde
    p=p_\delta + s + (\beta^{-1}(\widetilde x){-}1)(p_\delta +s).
  \end{align*}
  As $y$ and $s$ are small, we have that $\beta^{-1}(\widetilde x) \widetilde p \in
  \ove{B_{1/\epsi}}$ and we have proved that $\widetilde p =
  \beta(\widetilde x) (\beta^{-1}(\widetilde x) \widetilde p) \in
  C_{ \widetilde x\epsi}$. This implies that $(\widetilde x,\widetilde
  p)\in C_\epsi$, hence
  $(x,p_\delta)\in C^\circ_\epsi$.
    \item Case $d(x,\Omega)>\epsi^2$. We have that  
      $C_{\widetilde x\epsi}=\ove{B_{1/\epsi}}$ for all $\widetilde x$ close to $x$. For all
      $\delta>0$, choosing $p_\delta\in B_{1/\epsi}$ with
      $|p-p_\delta|\leq \delta$ we have $(x,p_\delta)\in C^\circ_\epsi$.
    \end{itemize}

As $\epsi <1/\sigma^*\leq 1/\sigma_*$, we can use \eqref{eq:H04} and deduce
for all $x \in \Rz^n$
that
\begin{align*}
  &B_{\sigma_*} = (\alpha(x)+\beta(x)) B_{\sigma_*}  \subset
    \alpha(x) B_{\sigma_*} + \beta(x)B_{1/\epsi} \\
  &\quad \subset
  C_{x\epsi} \subset  \alpha(x) B_{\sigma^*} +
  \beta(x)\ove{B_{1/\epsi}}\subset  (\alpha(x)+\beta(x))
  \ove{B_{1/\epsi}} = \ove{B_{1/\epsi}}.
\end{align*}
This proves that $\ove H_\epsi$ fulfills \eqref{eq:H04} in the whole
$\Rz^n$, but with $\sigma^*$ replaced by $ 1/\epsi $, 
namely,  \begin{equation}
  \label{eq:H04b}
  B_{\sigma_*}  \subset
  C_{x\epsi} \subset   \ove{B_{1/\epsi}}  \quad \forall x \in \Rz^n.
\end{equation}

\medskip

\noindent {\it Step 3: Solution to the approximating unconstrained problem}.
Since $ \ove H_\epsi$ fulfills  the assumptions of Proposition
\ref{prop:rn}  (with $\sigma^*$ replaced by
$1/\epsi$  in \eqref{eq:H04}),  the function 
\begin{align}
     v_\epsi(x) = \min \Bigg\{\int_0^{1} 
  \sigma_\epsi(\gamma(s),\gamma'(s)) \ {\d s}  \ \Big| \ \gamma \in  {\mathcal C}_x \Bigg\}\quad  \forall x
                  \in \Rz^n  \label{eq:vis5}
\end{align}
is the unique nonnegative viscosity solution to
\eqref{eq:1}--\eqref{eq:2} for $\Omega=\Rz^n$ and
$H$ replaced by $ \ove H_\epsi$. In
\eqref{eq:vis5}, the support function $\sigma_\epsi$ is defined as 
\begin{align*}
  \sigma_\epsi(x,q) &:= \sup\{q\cdot p \ | \  \ove H_\epsi(x,p)\leq
    0\}  \quad \forall
                      (x,q)\in \Rz^n \times \Rz^n.
\end{align*}
Owing to \eqref{eq:H04b} we have that
\begin{align}
  &\sigma_* |q| \leq \sigma_\epsi (x,q)\leq \frac{1}{\epsi}|q |\quad  \forall
    (x,q)\in \Rz^n \times \Rz^n,\label{eq:sigmae1}\\
  &\sigma_\epsi (x,q)=\sigma(x,q)\leq \sigma^*|q |\quad  \forall
    (x,q)\in \ove \Omega\times \Rz^n.\label{eq:sigmae2}
\end{align}
Moreover,  as
$C_{ \ove\pi (x)}=(\alpha(x)+\beta(x))C_{ \ove\pi (x)}\subset \alpha(x)C_{ \ove\pi (x)}+
\beta(x)\ove{B_{1/\epsi}}=C_{x\epsi}$ for all $x\in \Rz^n$ such that
$d(x,\Omega)\leq \epsi^2$, it follows that
  \begin{equation}
    \label{eq:proj}
    \sigma( \ove\pi (x),q) \leq
\sigma_\epsi(x,q)\quad \forall (x,q)\in \Rz^n\times \Rz^n \ \
\text{with} \ \ d(x,\Omega)\leq \epsi^2.
\end{equation} 

\medskip
\noindent {\it Step 4: Estimate and limit for $\epsi \to 0$}.
  For all $x,\, y \in \ove \Omega$, let $\gamma_y\in   {\mathcal
    C}_y $
  be optimal for $v_\epsi$ at $y$ and $\gamma \in W^{1,\infty}(0,1;\ove \Omega)$
  be a geodesic from $y$ to $x$. By
  considering the concatenated curve $\widetilde \gamma\in  {\mathcal C}_x $ given by
  $$\widetilde \gamma(s)=
\left\{
  \begin{array}{ll}
    \gamma_y(2s)\quad&\text{for}  \ s \in [0,1/2],\\[1mm]
    \gamma(2s-1)&\text{for}  \ s \in (1/2,1],
  \end{array}
  \right.
  $$ 
  we can bound
  \begin{align}
   & v_\epsi(x)-v_\epsi(y) \leq \int_0^{1}
    \sigma_\epsi(\widetilde\gamma(s),\widetilde \gamma'(s))\, \d s - \int_0^{1} \sigma_\epsi
    ( \gamma_y(s),\gamma'_y(s))\, \d s \nonumber\\
    &\quad = \int_{1/2}^{1}\sigma_\epsi
    (\widetilde \gamma(s),\widetilde \gamma'(s))\, \d s \stackrel{ \gamma(\cdot) \in
      \ove \Omega}{=} \int_{1/2}^{1}\sigma
    (\widetilde \gamma(s),\widetilde \gamma'(s))\, \d s
      \stackrel{\eqref{eq:sigmae2}}{\leq} \sigma^* \int_{0}^1 | \gamma'(s)|\, \d s\nonumber\\
    &\quad= \sigma^*{d_{\ove \Omega}(x,y)}
    \stackrel{\eqref{eq:E2}}{\leq}  \sigma^*L |x-y|.\label{eq:dopo}
  \end{align}
  By exchanging the roles of $x$ and $y$, this proves that
  $|\nabla v_\epsi|\leq \sigma^* L$ almost everywhere in
  $\Omega$.  As $v_\epsi (x) \leq \sigma^*L|x-x_0|$  for some
  fixed $x_0\in V_0$ and  for all $x \in
  \ove \Omega$, a  diagonal extraction argument allows   one  to find $v\in
  C(\ove \Omega)$ and a not relabeled subsequence $v_\epsi$ such that 
  $$v_\epsi \to v \quad \text{locally uniformly in} \ \ove \Omega.$$

It is a standard matter to prove that  $v$ is a viscosity
solution to the constrained problem
\eqref{eq:1}--\eqref{eq:2}. Indeed, for all $x\in \Omega\setminus \ove{V_0}$ and all $\varphi \in
C^1(\Rz^n)$ with $\varphi (x) - v(x)= \min (\varphi - v)$
 locally  we use
\cite[Prop.~4.3]{userguide} in order to find $\Omega\setminus \ove{V_0}\ni
x_{\epsi}\to x$ and $\varphi_{\epsi}\in C^1(\Rz^n)$ with $\varphi_\epsi(x_\epsi) -
v_\epsi(x_\epsi)= \min (\varphi_\epsi - v_\epsi)$ 
locally  and $\nabla \varphi_{ \epsi} (x_\epsi) \to
\nabla \varphi (x)$. As $ \ove H(z,\cdot)= \ove H_\epsi(z,\cdot)$ for
 all  $z \in \Omega
\setminus \ove{V_0}$  and $\ove H$ is continuous, see
\eqref{eq:Hcontinuous}, 
one gets $0 \geq \lim_{\epsi \to 0} \ove H_\epsi (x_\epsi,\nabla
\varphi_{ \epsi}
(x_\epsi)) =    \lim_{\epsi \to 0} \ove H (x_\epsi,\nabla
\varphi_{ \epsi}
(x_\epsi)) = \ove H (x,\nabla \varphi
(x)) $, proving that $v$ is a viscosity subsolution to
\eqref{eq:1} in $\Omega\setminus \ove{ V_0}$. Similarly one proves that $v$ is a viscosity
supersolution, hence a viscosity solution.

\medskip
\noindent {\it Step 5: Representation formula.} We now prove that the
limit $v \in C(\ove \Omega)$ is represented as in
\eqref{eq:vis4}. Fix $x \in \ove \Omega$.  If $x \in \ove
V_0$ the formula trivially holds,  hence  we can assume $x \in
\ove \Omega \setminus \ove {V_0}$. Let   $\gamma \in
 {\mathcal C}_x $ with $\gamma(s) \in \ove \Omega$ for all $s \in [0,1]$. Then,
\begin{equation}
  \label{eq:liminf}
  v(x) = \lim_{\epsi \to 0}v_\epsi(x) \leq \int_0^{1}\sigma_\epsi(\gamma(s),\gamma'(s))\, \d
  s \stackrel{\gamma(\cdot)\in \ove \Omega}{=} \int_0^{1}\sigma(\gamma(s),\gamma'(s))\, \d
  s
\end{equation}
so that the right-hand side of \eqref{eq:vis4} is an upper bound for
$v(x)$.

We are left with checking that the minimum in \eqref{eq:vis4} is
actually attained, namely, that there exists an optimal curve  in
$\ove \Omega $  realizing
the equality in \eqref{eq:liminf}. To this aim, let
$\gamma_{x\epsi} \in  {\mathcal C}_x $ be optimal curves for $x$
under $\sigma_\epsi$  and note that the corresponding lengths
$\ell_{\epsi}$ are uniformly bounded as
\begin{equation*}
  \ell_{\epsi}= \int_0^1 |\gamma_{x\epsi}'(s)|\, \d s
  \stackrel{\eqref{eq:sigmae1}}{\leq}\frac{1}{\sigma_*}  \int_0^1
  \sigma_\epsi(\gamma_{x\epsi}(s), \gamma_{x\epsi}'(s))\, \d s = \frac{v_\epsi(x)}{\sigma_*}
\end{equation*}
and $v_\epsi(x) \to v(x)$. In particular, $\ell_{\rm
  max}:=\sup_\epsi \ell_{\epsi}<\infty$. On the other hand, the
lengths $\ell_\epsi$ are bounded below by $d(x,\ove V_0)>0$.
We indicate by $\widetilde \gamma_{x\epsi} \in
W^{1,\infty}(0,\ell_{\epsi};\Rz^n)$ the 
arc-length reparametrizations of $\gamma_{x\epsi}$, see
\cite[Lemma~1.1.4(b), p.~25]{Ambrosio}, and trivially extend them to
$[0,\ell_{\rm \max}]$ without changing notation.  As the curves $\widetilde\gamma_{x\epsi}$ are
uniformly Lipschitz continuous and $\widetilde\gamma_{x\epsi}(\ell_{\rm max})=x$, one can extract not relabeled
subsequences in such a way that $\ell_{\epsi} \to \ell \in [d(x,\ove V_0),
\ell_{\rm max}]$ and $\widetilde \gamma_{x\epsi} \to \widetilde \gamma_x$ uniformly in $
[0,\ell_{\rm max}]$ and weakly* in $W^{1,\infty}(0,\ell_{\rm
  max};\Rz^n)$. Note that
$|\widetilde \gamma_x '|\leq 1$ a.e. in $(0,\ell )$ and $
\widetilde \gamma_x ' =0 $ a.e. in $(\ell ,\ell_{\rm
  max})$. Hence,  
$\widetilde \gamma_x(\ell )=\widetilde \gamma_x(\ell_{\rm max})=x$, and $\ove{V_0} \ni
\widetilde\gamma_{x\epsi}(0)\to \widetilde\gamma_x(0)$. As $\ove{V_0} $
is closed, the latter convergence ensures that $\widetilde\gamma_x(0)\in
\ove{V_0} $, as well. 

Let us now
check that $\widetilde\gamma_x(s)\in \ove \Omega$ for all $s\in [0,
\ell ]$. To
this aim, define
$$A_\epsi:=\{s \in  [0,\ell_{\epsi}]\ | \  d(\widetilde
\gamma_{x\epsi}(s),\ove \Omega)\geq \epsi^2\}.$$
For almost all  $s \in A_\epsi$ we have that
$$\sigma_\epsi (\widetilde
\gamma_{x\epsi}(s),\widetilde \gamma_{x\epsi}'(s))
=\frac{1}{\epsi}|\widetilde \gamma_{x\epsi}'(s)|=
\frac{1}{\epsi}.$$
This allows us to estimate $|A_\epsi|$ as follows
$$ \frac{1}{\epsi}|A_\epsi| = \int_{A_\epsi} \sigma_\epsi (\widetilde
\gamma_{x\epsi}(s), \widetilde
\gamma_{x\epsi}'(s))\, \d s \leq \int_0^{\ell_\epsi} \sigma_\epsi (\widetilde
\gamma_{x\epsi}(s), \widetilde
\gamma_{x\epsi}'(s))\, \d s = v_\epsi(x) \to v(x),$$
implying that
$$
|A_\epsi|\leq c\epsi
$$
for some $c>0$ depending on  
$v(x)$ only. Hence, the length of the portion of the
arc-length parametrized curve $\widetilde
\gamma_{x\epsi}$  lying   in the complement of $ \ove \Omega+ {B_{\epsi^2}}$ is smaller than $
c\epsi$. This in particular entails that
$$  d(\widetilde \gamma_{x\epsi}(s),  \ove \Omega+ {B_{\epsi^2}})\leq
 c\epsi  \quad \forall s \in[0,\ell_{\epsi}]$$
and by the triangle inequality we get
$$  d(\widetilde \gamma_{x\epsi}(s),\ove \Omega)\leq  d(\widetilde \gamma_{x\epsi}(s),\ove
  \Omega+ {B_{\epsi^2}}) +\epsi^2 \leq    c
  \epsi +\epsi^2  \quad \forall s \in[0,\ell_{\epsi}].
  $$
By passing to the limit as $\epsi \to 0$ we obtain that $\widetilde \gamma_x(s) \in
\ove \Omega$ for all $s \in
 [0,\ell ]$.

Let now $\delta\in (0,\ell )$ be given and choose  elements
$\epsi_j$ of the above extracted sequence such that $\sum_j \epsi_j
<\delta$.  For $j$ large enough, one has that $c\epsi_j+\epsi^2_j
<\ove \epsi$. Define $A :=\cup_{j}A_{\epsi_j}$ so that 
$|A |\leq \sum_j |A_{\epsi_j}| \leq c\delta $. By using \eqref{eq:proj} we find that 
\begin{align*}
v(x)&=\lim_{j \to \infty}
                 v_{\epsi_j}(x)=\lim_{j \to
           \infty}\int_0^{\ell_{\epsi_j}}\sigma_{\epsi_j}(\widetilde\gamma_{x\epsi_j}(s),\widetilde\gamma_{x\epsi_j}'(s))\,
           \d s  \\
  & \geq \liminf_{j \to \infty}  \int_{[0, \ell_{\epsi_j}]\setminus A_{\epsi_j}}\sigma_{\epsi_j}(\widetilde\gamma_{x\epsi_j}(s),\widetilde\gamma_{x\epsi_j}'(s))\,
    \d s  \\
  &\hspace{-2.2mm}\stackrel{ \eqref{eq:proj}}{\geq}
     \liminf_{j \to \infty} \int_{[0, \ell_{\epsi_j} ]\setminus A_{\epsi_j} }\sigma( \ove\pi (\widetilde\gamma_{x\epsi_j}(s)),\widetilde\gamma_{x\epsi_j}'(s))\, \d
    s  \\
  &\geq \liminf_{j \to \infty} \int_{[0, \ell -\delta]\setminus A }\sigma( \ove\pi (\widetilde\gamma_{x\epsi_j}(s)),\widetilde\gamma_{x\epsi_j}'(s))\, \d
                            s \\
  & \geq \int_{[0, \ell-\delta ]\setminus A }\sigma( \ove\pi (\widetilde\gamma_x(s)), \widetilde\gamma_x'(s))\, \d
    s   \\
  &\geq \int_0^\ell\sigma( \ove\pi (\widetilde\gamma_x(s)), \widetilde\gamma_x'(s))\, \d
     s    -\sigma^* \delta - \sigma^*|A| \\
  & \hspace{-3.6mm}\stackrel{\widetilde\gamma_x(\cdot) \in \ove \Omega}{\geq}
    \int_0^\ell\sigma(\widetilde\gamma_x(s), \widetilde\gamma_x'(s))\, \d
     s - \sigma^*(1+c) \delta
\end{align*}
where we have used the continuity of $\sigma$ and $ \ove\pi $, the
convexity
of $\sigma(x,\cdot)$, see Lemma~\ref{lemma:sigma}, and the classical
\cite[Thm.~3.23, p.~96]{Dacorogna}. As $\delta$ is arbitrary, by defining
$\gamma_x\in  {\mathcal C}_x $ as $\gamma_x(s) = \widetilde
\gamma_x(\ell s)$ for all $s\in [0,1]$, this proves that 
$$v(x)\geq \int_0^\ell\sigma(\widetilde\gamma_x(s), \widetilde\gamma_x '(s))\, \d
     s = \int_0^1 \sigma(\gamma_x(s),\gamma_x'(s))\, \d
     s.$$
Together with the upper bound
\eqref{eq:liminf}, we have checked that $\gamma_x$ is optimal for $x$ and the
representation formula \eqref{eq:vis4} holds.

\medskip
\noindent {\it Step 6: Boundary behavior.}
Assume by contradiction that
$v $ is not a constrained viscosity solution. In particular,
 assume that  there exist $x \in \partial \Omega \setminus
\ove{V_0}$  and  $\varphi
\in C^1(\Rz^n)$ such that $\varphi - v$ has a local maximum at
$x$ but $H(x,\nabla \varphi(x))<0$. Let $\gamma \in  {\mathcal C}_x $
realize the minimum in \eqref{eq:vis4} and assume with no loss of
generality that $|\gamma'|=  \omega   >0$ a.e. For all $r>0$ we can find
$t\in (0,1)$ such that $\gamma(s)\in B_r(x)$ for all $s \in
[t,1]$. Owing to the continuity of $H$ and $\nabla \varphi$ and
possibly choosing a larger $t$, we find $\delta>0$ such that $H(\gamma(s),
\nabla \varphi(\gamma(s)) ) <-\delta$ for all $s \in [t,1]$. In particular, there exists $\lambda \in (0,1)$ with $H(\gamma(s), \nabla
\varphi(\gamma(s))/\lambda )\leq 0$ for all $s \in [t,1]$. Using the very definition of
$\sigma$ this entails that
\begin{equation}\label{eq:lambda}
  \gamma'(s)\cdot
\nabla \varphi(\gamma(s)) \leq \lambda \sigma
(\gamma(s),\gamma'(s))\quad \text{for a.e.} \ s \in (t,1).
\end{equation}
We can now compute
\begin{align}
  &\varphi(x)- \varphi(\gamma(t)) =  \int_t^1\gamma'(s)\cdot \nabla
    \varphi(\gamma(s))\, \d s  \stackrel{\eqref{eq:lambda}}{\leq}   \lambda \int_t^1\sigma
    (\gamma(s),\gamma'(s))\, \d s  \nonumber\\
  &\quad =  \int_t^1\sigma
    (\gamma(s),\gamma'(s))\, \d s + (\lambda -1)  \int_t^1\sigma
    (\gamma(s),\gamma'(s))\, \d s   \nonumber\\
  &\quad \leq \int_0^1\sigma
    (\gamma(s),\gamma'(s))\, \d s  - \int_0^t\sigma
    (\gamma(s),\gamma'(s))\, \d s +(\lambda -1)(1-t)  \omega  \sigma_*  \nonumber\\[2mm]
  &\quad \leq v(x) - v(\gamma(t)) +(\lambda
    -1)(1-t)  \omega  \sigma_*< v(x) - v(\gamma(t)) .\label{eq:can}
\end{align}
This implies that $\varphi(\gamma(t)) - v(\gamma(t)) > \varphi(x)
- v(x)$. As $r>0$ is arbitrary, we  hence  proved that $x$ is not a local maximum
of $\varphi - v$, which is a contradiction.
\end{proof}

\begin{remark}[Alternative assumptions, continued]\label{rem:seconda}
   In case $\Omega $ is convex the proof of
  Proposition~\ref{prop:representation} can be simplified by skipping Step 1 and
  arguing directly on the extension $H_\epsi(x,p)=H( \ove\pi (x),
  b_\epsi(x)p)$ where $ \ove\pi $ is again  the projection on $\ove \Omega$ and
  $b_\epsi(x)=a_\epsi(d(x,\ove \Omega)/\epsi^2)$ for some
  $a_\epsi: \Rz \to [\epsi, 1]$ smooth with $a_\epsi(r) =1 $
  for $r \leq 1$ and $a_\epsi(r)=\epsi$ for $r\geq 2$. In fact,
  in this case one only needs $ \ove\pi $ to be a
  continuous map from $\Rz^n$ to $\ove \Omega$ coinciding with the
  identity on $\ove \Omega$. Hence, one could also consider some
  nonconvex $\ove \Omega=h(\ove D)$ as indicated in
  Remark~\ref{rem:prima} by replacing $ \ove\pi $ in $H$
  by $h\circ  \ove\pi _{D} \circ h^{-1}$,  where now $ \ove\pi _D$ is the
  projection on $ \ove D $. 
\end{remark}

Let us record some basic estimates from the proof of
Proposition~\ref{prop:representation}   in the following.

\begin{proposition}[Estimates on $v$]\label{prop:estimates}
  Under the assumptions   of 
     Proposition~\ref{prop:representation},  let  $x_0 \in V_0$ and
   $v$ be given by \eqref{eq:vis4}. Then, the length $\ell_x$ of any
   optimal curve $\gamma_x\in  {\mathcal C}_x $ with $\gamma_x(\cdot) \in \ove
   \Omega$ fulfills
   \begin{equation}
     \label{eq:est0}
     \ell_x \leq \frac{v(x)}{\sigma_*} \quad \forall x \in \ove \Omega.
   \end{equation}
   Moreover, we have that 
  \begin{align}
    &|\nabla v| \leq \sigma^*L\quad \text{a.e. in} \ \Omega, \label{eq:est1}\\
    &\forall R>0: \ v(x) \leq c_{ 1 }(R) \quad \forall x \in \ove \Omega \cap
      B_R, \label{eq:est2}      
  \end{align}
  where $c_{ 1}(R):=\sigma^* L (|x_0|+R) $ for all $R>0$.
\end{proposition}

\begin{proof}
   The bound \eqref{eq:est0} on the length $\ell_x$ follows by 
   \eqref{eq:sigma}  as
   $$v(x) = \int_0^{ 1 }\sigma(\gamma_x(s), \gamma_x'(s))\, \d s \geq
   \sigma_* \ell_{ x} .$$
  The Lipschitz bound \eqref{eq:est1}  can be obtained  by repeating the
  argument leading to \eqref{eq:dopo}  for the function $v$,
  using \eqref{eq:vis4}. 

 For \eqref{eq:est2}, fix any $x \in \ove \Omega \cap B_R$ and let
 $ \haz \gamma  \in  {\mathcal C}_x $ be a geodesic in $\ove
 \Omega$ from $x_0$ to $x$. We have
 \begin{align*}
   v(x)&\stackrel{\eqref{eq:vis4}}{\leq}
   \int_0^1\sigma( \haz \gamma  (s),  \haz \gamma'  (s))\, \d s \stackrel{\eqref{eq:sigma}}{\leq}
   \sigma^*\int_0^1| \haz \gamma'  (s)|\, \d s \\
   &\hspace{1.6mm} =\sigma^*d_{\ove \Omega}(x_0,x)
   \stackrel{\eqref{eq:E2}}{\leq} \sigma^* L |x-x_0|\leq \sigma^* L (|x_0|+R).\qedhere
 \end{align*} 
\end{proof}

In view of the coupling with the elliptic problem for the activation
field $u$ in Section \ref{sec:elastic}, the regularity of the sublevels
$V(t)$ plays a crucial role. We have the following.
  
  \begin{proposition}[John regularity of $V(t)$]\label{prop:john}
    Under the assumptions of
    {Proposition~\ref{prop:representation}}, assume additionally
    \eqref{eq:H00u} and that  
 $v$ is given by \eqref{eq:vis4}.  Then, for all $T>0$ there
 exists $\kappa(T)\in (0,1]$ such that 
  $V(t) = \{x
  \in  \Omega\ | \ v(x)<t\}$ are  John domains with respect to
  $x_0$ with John constant $\kappa(T) $ for all
  $t\in (0,T]$, where  $\kappa(T) $ depends  on  $\sigma_*$,
 $\sigma^*$, $\ove \epsi$,   $L$,  $\|\nabla \eta \|_{L^\infty}$, $\kappa_0$,
 $d(x_0,\partial \Omega)$,  $d(x_0,\partial V_0)$, ${\rm
   diam}(V_0)$  (although not indicated
 in the notation), and $T$, and $\kappa(T) \to 0$ as $T\to
 \infty$. Moreover, we have that
   \begin{align} 
    &d_H(V(t),V(s))\leq \frac{1}{\sigma_*}(t-s) \quad
      \forall 0<s\leq t   \label{eq:lipschitz}.
   \end{align}
   Finally, for all $T>0$  there exists $\mu(T)\in (0,1)$ and $c_H(T)>0$ such that
   \begin{align}  
    \label{eq:hoelder}
  &  |V(t)\setminus V(s)|\leq c_H(T) (t-s)^{\mu(T)} \quad \forall 0<s\leq t \leq
    T,
   \end{align}
   where $\mu(T)$ and $c_H(T)$ have the same dependencies  as  $\kappa(T)$.
\end{proposition}

In case $\Omega=\Rz^n$, the John regularity of $V(t)$ is proved in
\cite[Thm.~1.1]{john}. In  the  unconstrained case,  the  assumptions on $H$
can be weakened and the constant $\kappa$ turns out to  depend 
only on $\kappa_0$,  $\sigma_*$, and $\sigma^*$  and, in particular, is independent of $T$. The constrained case $\Omega \not = \Rz^n$ is more
delicate, as the presence and the regularity of the boundary $\partial \Omega $ plays a
role. Note nonetheless that John regularity is sharp for problem
\eqref{eq:1}--\eqref{eq:2}, see \cite{john} for counterexamples.

\begin{proof}[Proof of Proposition~\ref{prop:john}]
We divide the proof into steps. At first, we find an arc-length
parametrized curve $\gamma$ in
$\ove \Omega$ such that $d(\gamma(s),\partial V(t)\setminus \partial
\Omega) \geq \ove \kappa s$ for some $\ove \kappa>0$ (Step 1, 
Figure~\ref{fig:2}).  Then,
starting from $\gamma$ we find a second curve $\widetilde \gamma$ such
that $d(\widetilde \gamma(s),\partial V(t)) \geq \widetilde \kappa s$
for some $0<\widetilde \kappa\leq \ove \kappa$ (Step 2,
Figure~\ref{fig:3}).   Finally,  we reparametrize
$\widetilde \gamma$ by arc  length and conclude that $V(t)$ is a  John
domain (Step 3). The Lipschitz continuity of  the set-valued map
 $t \mapsto V(t)$
with respect to the Hausdorff distance and the H\"older continuity
with respect to the $L^1$ distance are proved in Step 4 and Step 5,
respectively.

\medskip

\noindent {\it Step 1: The curve $\gamma$.} For all $T>0$ fixed, let $t \in (0,T]$ and $x \in
V(t)$. Let $\haz \gamma_1\in  {\mathcal C}_x $ with $\haz \gamma_1(\cdot) \in \ove \Omega$ be
optimal for $x$, indicate by $ \gamma_1:[0,\ell_{\gamma_1}]\to
\ove \Omega$ its arc-length reparametrization (see, e.g.,
\cite[Lemma~1.1.4(b), p.~25]{Ambrosio}), and set $y =\gamma_1(0)\in \ove
{V_0}$. Moreover, using the fact that $\ove{V_0}$ is a John domain, set $\gamma_2:[0,\ell_{\gamma_2}]\to
\ove{V_0}$ to be such that $\gamma_2(0)=y$,
$\gamma_2(\ell_{\gamma_2})=x_0$, $|\gamma_2'|=1$ a.e., and
$d(\gamma_2(s),\partial V_0)\geq \kappa_0s$ for all $s\in
[0,\ell_{\gamma_2}]$. Set $\ell_\gamma=\ell_{\gamma_1}+
\ell_{\gamma_2}$ and define $\gamma:[0,\ell_\gamma]\to \ove \Omega$ by
$$\gamma(s):=
\left\{
  \begin{array}{ll}
    \gamma_1(\ell_{\gamma_1}-s) \quad&\text{for} \ s \in
                                       [0,\ell_{\gamma_1}],\\
    \gamma_2(s-\ell_{\gamma_1}) \quad&\text{for} \ s \in
                                       (\ell_{\gamma_1},\ell_\gamma],
  \end{array}
  \right.
  $$
  see Figure~\ref{fig:2}. 
  In particular, $\gamma(0)=x$, $\gamma(\ell_\gamma)=x_0$, and
  $|\gamma'|=1$ a.e.  Note that
   \begin{equation}\label{eq:lmax}\ell_\gamma \leq  \ell(T):= 
 \frac{d(x_0,\partial V_0)}{\kappa_0} +
 \frac{T}{\sigma_*}.
\end{equation}

We devote the rest of this step to proving that,
  setting
  \begin{equation}
    \ove \kappa =\frac{\sigma_*}{ 2L\sigma^*+ \sigma_*}   
  \kappa_0,\label{eq:CJohn}
  \end{equation}
  one has that
 \begin{equation}
  d(\gamma(s),\partial V(t)\setminus \partial \Omega)\geq
 \ove \kappa s \quad \forall s \in [0,\ell_{\gamma}].  \label{eq:elisa0}
\end{equation}
This follows from adapting 
 the argument of
 \cite[Thm.~1.1]{john}.  We argue by proving \eqref{eq:elisa0}
 separately on the three intervals $[0,\ell_{\gamma_1}]$,
 $[\ell_{\gamma_1},s^*]$, and $[s^*,\ell_\gamma]$, where $s^* $ is
 defined in \eqref{eq:sstar} below. 

At  first, we check that 
  \begin{equation}
  d(\gamma(s),\partial V(t)\setminus \partial \Omega)\geq
  \frac{\sigma_*}{ L \sigma^*}s \quad \forall s \in [0,\ell_{\gamma_1}].  \label{eq:elisa1}
\end{equation}
To this end, fix $s\in [0,\ell_{\gamma_1}]$ and  let
$\rho:[0,\ell_\rho]\to  \ove \Omega $ be any arc-length parametrized curve connecting $
\gamma(s)$ to $\partial V(t)\setminus \partial \Omega$, namely, such
that  $\rho(0) \in \partial V(t) \setminus \partial \Omega$ and
$\rho(\ell_\rho)=\gamma(s)$. The curve resulting from following 
$\gamma_1(r)$ for $r \in [0,\ell_{\gamma_1} -s ]$ and then
$\rho(\ell_{\gamma_1}+ \ell_\rho -t)$ for $t \in (\ell_{\gamma_1}-s,
\ell_{\gamma_1} + \ell_\rho -s]$ connects $\ove{V_0}$ to $\partial
V(t)$. Hence, formula \eqref{eq:vis4} and the fact that $v(\cdot ) =t$
on $\partial V(t)\setminus \partial \Omega$ give
\begin{align}
  t \leq \int_0^{\ell_{\gamma_1}-s}\sigma(\gamma_1(r),\gamma_1'(r))\,
  \d r + \int_{\ell_{\gamma_1}-s}^{\ell_{\gamma_1} + \ell_\rho -s}
  \sigma(\rho(\ell_{\gamma_1} { +} \ell_\rho {-}r), -\rho'(\ell_{\gamma_1} {+}
  \ell_\rho {-}r)) \, \d r .\label{eq:inequa1}
\end{align}
On the other hand, we have that
\begin{align}
  &\int_0^{\ell_{\gamma_1} -s} \sigma(\gamma_1(r), \gamma_1'(r)) \, \d
    r + \int_{\ell_{\gamma_1} -s}^{\ell_{\gamma_1}} \sigma(\gamma_1(r), \gamma_1'(r)) \, \d
    r \nonumber\\
  &\quad = \int_0^{\ell_{\gamma_1}} \sigma(\gamma_1(r), \gamma_1'(r)) \, \d
    r = v(x) <t .\label{eq:inequa2}
\end{align}
Putting the inequalities \eqref{eq:inequa1}--\eqref{eq:inequa2}
together we obtain
$$\int_{\ell_{\gamma_1} -s}^{\ell_{\gamma_1}} \sigma(\gamma_1(r), \gamma_1'(r)) \, \d
    r < \int_{\ell_{\gamma_1}-s}^{\ell_{\gamma_1} + \ell_\rho -s}
  \sigma(\rho(\ell_{\gamma_1} { +} \ell_\rho {-}r), -\rho'(\ell_{\gamma_1} {+}
  \ell_\rho {-}r)) \, \d r .$$ 
By using the bounds \eqref{eq:sigma} on 
$ \sigma$, the latter gives
$$ \sigma_*{s} <\sigma^* \ell_\rho ,$$
implying that the length $\ell_\rho$ of any curve  in $\ove
\Omega$  connecting $
\gamma(s) $ to $\partial V(t)\setminus \partial \Omega$ is strictly bounded  from  below
by $\sigma_* s/\sigma^*$. Taking $\rho$ to be a  geodesic   in
$\ove \Omega  $  from
$\gamma(s)$ to $x \in \partial V(t)\setminus \partial \Omega$ so
that $|\gamma(s) -x|=d(\gamma(s),\partial V(t)\setminus \partial \Omega)$, this proves
that
$$ d(\gamma(s), \partial V(t)\setminus \partial \Omega)=
|\gamma(s) - x| \geq \frac{1}{L} d_{\ove \Omega}(\gamma(s),x)=
\frac{\ell_\rho}{L}\geq \frac{\sigma_*}{L\sigma^*}s,$$ 
giving  the lower
bound \eqref{eq:elisa1}.  As $\ove \kappa \leq \sigma_*/\sigma^*$ (recall
that $\kappa_0 \leq 1$), this implies  the lower bound \eqref{eq:elisa0}
 for $s\in
[0,\ell_{\gamma_1}]$.
\begin{figure}
  \centering
  \pgfdeclareimage[width=95mm]{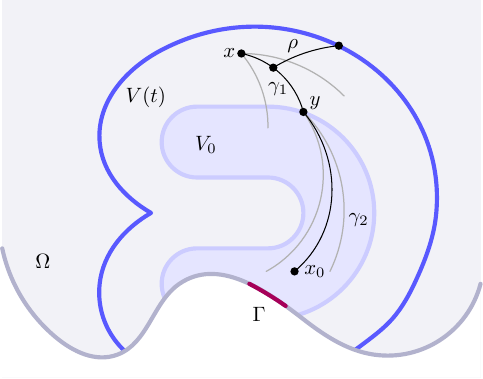}{figure2}
  \pgfuseimage{figure2}
  \caption{The curve $\gamma$ is obtained by concatenating $\gamma_1$
    and $\gamma_2$. The twisted cone around $\gamma_1$ is described by
  \eqref{eq:elisa1}, whereas the one around $\gamma_2$ follows from $V_0$
being a John domain.}
  \label{fig:2}
\end{figure}

Consider now the case $s\in [\ell_{\gamma_1}, s^*]$ with
\begin{equation}
  s^*:=\min\left\{\ell_{\gamma_1} \left(1+
\frac{\sigma_*}{ 2L\sigma^* }\right),
\ell_{\gamma}\right\}.\label{eq:sstar}
\end{equation}
In particular,  $s-\ell_{\gamma_1} \leq s^*-\ell_{\gamma_1} \leq
 \sigma_*\ell_{\gamma_1}/( 2L \sigma^*  )
$. This implies that 
\begin{align}
   |y - \gamma (s)|&= |y -  \gamma_2(s-\ell_{\gamma_1} )|
   =|\gamma_2(0) -  \gamma_2(s-\ell_{\gamma_1} )|\leq s -
  \ell_{\gamma_1} \leq  {\sigma_*\ell_{\gamma_1}}/{( 2L \sigma^*  )}, \nonumber
  \end{align}
    where we also used that $|\gamma_2(0) -
    \gamma_2(s-\ell_{\gamma_1} )|\leq s - \ell_{\gamma_1} $, as
    $\gamma_2$ is arc-length parametrized. As $d(y , \partial
    V(t)\setminus \partial \Omega) \geq
    {\sigma_*\ell_{\gamma_1}}/{( L  \sigma^*)} $ from
    \eqref{eq:elisa1},  the triangle inequality allows us to  conclude that
    \begin{align*}
      \frac{\sigma_*\ell_{\gamma_1}}{ L \sigma^* }  &\leq
                                                              d(y,
                                                              \partial
                                                              V(t) 
      \setminus \partial \Omega)\leq  d(\gamma(s), \partial V (t) 
        \setminus \partial \Omega) + |y-\gamma(s)| \\
      & \leq  d(\gamma(s), \partial V  (t) 
      \setminus \partial \Omega) + \frac{\sigma_*\ell_{\gamma_1}}{
        2L \sigma^*  }.
    \end{align*}
   This implies that 
    $d(\gamma (s) , \partial
    V(t)\setminus \partial \Omega) \geq
    {\sigma_*\ell_{\gamma_1}}/{( 2L \sigma^*  )} $ and \eqref{eq:elisa0}
    holds for $s \in [\ell_{\gamma_1},s^*]$, as well  as
    \begin{align*}
      &d(\gamma (s) , \partial
    V(t)\setminus \partial \Omega)
        \geq\frac{\sigma_*\ell_{\gamma_1}}{ 2L \sigma^* }  =
\frac{\sigma_*}{ 2L\sigma^*+ \sigma_* }\ell_{\gamma_1}\left(
  1+\frac{\sigma_*}{ 2L \sigma^* }\right)
\stackrel{\eqref{eq:CJohn}}{\geq }  \ove \kappa s^*\geq   \ove \kappa s. 
    \end{align*}
    
We are hence left with considering the case $s\in [s^*,
\ell_{\gamma}]$. Note that this case is only relevant if 
$s^*<\ell_\gamma$, hence for $s^*=\ell_{\gamma_1
}(1+\sigma_*/( 2L \sigma^* ))$.  As $s > \ell_{\gamma_1}$, we have
that $\gamma(s) \in V_0$.  We use that $\ove{V_0} $ is a John
domain with John center $x_0$ and John constant $\kappa_0$ to obtain
$$d(\gamma(s), \partial V(t)\setminus \partial \Omega) \geq  d(\gamma(s),
\partial V_0)=d(\gamma_2(s-\ell_{\gamma_1}), \partial V_0) \geq
\kappa_0 (s- \ell_{\gamma_1}).$$
The lower bound \eqref{eq:elisa0}  for $s \in [s^*,\ell_\gamma]$
follows then from    
\begin{align}
  &\kappa_0 (s- \ell_{\gamma_1}) =  \kappa_0 \left( 1-
  \frac{\ell_{\gamma_1}}{s}\right)s \geq \kappa_0\left( 1-
    \frac{\ell_{\gamma_1}}{s^*}\right)s  
  = \kappa_0\frac{\sigma_*}{ 2L \sigma^*+ \sigma_*  }s
    \stackrel{\eqref{eq:CJohn}}{\geq }   \ove \kappa s.\nonumber
\end{align} 
\medskip

 \noindent {\it Step 2:  The curve  $\widetilde \gamma$.}

   To prove that $V(t)$ is a John domain, we need to complement  
  \eqref{eq:elisa0} by providing a lower bound for the distance of
  $\gamma( s)$ from $\partial
  V(t) \cap \partial  \Omega$. This however may require  modifying
   the
  curve $\gamma$, as a priori $\gamma$ may get close and even touch $\partial V(t)$ at points in $\partial
  V(t) \cap \partial  \Omega$,  see Figure~\ref{fig:3}. We hence
  define a new curve $\widetilde \gamma$ by displacing $\gamma$ in the
  inward normal direction in a neighborhood of the boundary $\partial
  V(t) \cap \partial  \Omega$.  This may also   require  reducing  the John constant.
  \begin{figure}
  \centering
  \pgfdeclareimage[width=95mm]{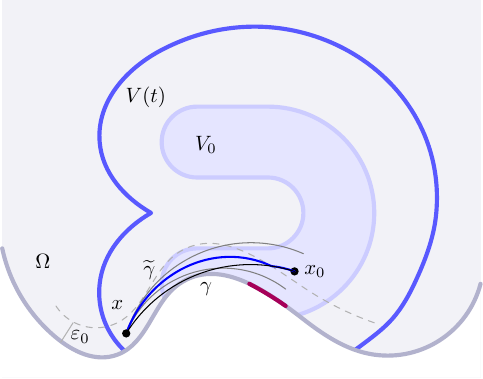}{figure3}
  \pgfuseimage{figure3}
  \caption{Construction of the curve $\widetilde \gamma$ (Step
    2). The curve $\gamma$ touches the boundary $\partial
    \Omega$. Within the $\varepsilon_0$-neighborhood of $\partial
    \Omega$ (indicated by the dashed curve), the curve
    $\widetilde\gamma$ is obtained by displacing $\gamma$ in the
    inward normal direction. This allows us to fit a twisted cone around
    $\widetilde \gamma$ contained in $V(t)$.}
  \label{fig:3}
\end{figure}

  Let us start by recalling that  the projection  $  \pi:(\partial \Omega+ B_{\ove \epsi })\cap \ove \Omega \to \partial
  \Omega  $ and the normal
outward-pointing vector field $\nu: \partial \Omega \to {\mathbb S}^{n-1}$ 
are well-defined and Lipschitz continuous \cite[Cor.~2]{Lewicka} and
set $\eta:=-\nu\circ \pi$, which is Lipschitz continuous, as well.

  Fix $0<\epsi_0 < \min\{\ove \epsi/2,d(x_0,\partial
  \Omega)\}$ and set $\widetilde \kappa>0$ so small that
  \begin{equation}
    \widetilde \kappa <\min\left\{\frac{\ove
      \kappa}{2},\frac{\epsi_0}{   \ell(T)  }, \frac{1}{ 2(\|
      \nabla \eta\|_\infty    \ell(T)  +    \ell(T)   /\epsi_0 +1)
    }\right\}\label{eq:tildek}
\end{equation}
 for $\ell(T)$ defined in \eqref{eq:lmax}.  
We define $\widetilde \gamma : [0,\ell_\gamma]\to \Omega$ as 
  $$\widetilde \gamma (s):= \gamma(s) + 
  \eta(\gamma(s))  \big(\epsi_0 - d(\gamma(s) ,\partial \Omega)\big)^+
  \frac{\widetilde \kappa}{\epsi_0} s.$$
  Note that   $\widetilde \gamma(0)=
  \gamma(0)=x$, and, using also $\gamma(\ell_\gamma)=x_0$ and the fact
  that $\epsi_0<d(x_0,\partial \Omega)$, we have that 
  $$\widetilde \gamma(\ell_\gamma) = x_0+  
  \eta(x_0)  \big(\epsi_0 - d(x_0,\partial \Omega)\big)^+
  \frac{\widetilde \kappa}{\epsi_0} \ell_\gamma= x_0.$$ Moreover, 
  $\widetilde \gamma$
  is Lipschitz continuous. 
  In fact, for all $s,\,t \in (0,\ell_\gamma)$
  we have that
  \begin{align*}
    &\widetilde \gamma(t) - \widetilde \gamma (s) = \gamma(t) - \gamma
      (s)\\
    &\quad +\big(\eta(\gamma(t))-\eta(\gamma(s))\big)  \big(\epsi_0 - d(\gamma(t) ,\partial \Omega)\big)^+
      \frac{\widetilde \kappa}{\epsi_0} t \\
    &\quad +
    \eta(\gamma(s))\Big(\big(\epsi_0 - d(\gamma(t) ,\partial
    \Omega)\big)^+-\big(\epsi_0 - d(\gamma(s) ,\partial \Omega)\big)^+ \Big)
      \frac{\widetilde \kappa}{\epsi_0} t\\
    &\quad +  \eta(\gamma(s)) \big(\epsi_0 - d(\gamma(s) ,\partial \Omega)\big)^+ \frac{\widetilde \kappa}{\epsi_0} (t-s).
  \end{align*}
  This ensures that
\begin{align*}
    &|\widetilde \gamma(t) - \widetilde \gamma (s) |\leq |\gamma(t) -
      \gamma(s)|+ \widetilde \kappa \left(  \|\nabla \eta\|_{L^\infty}       \ell(T)  +   \ell(T)  /\epsi_0\right) |\gamma(t) - \gamma(s)| + \widetilde
      \kappa |t-s|
  \end{align*}
and, in particular,
  \begin{equation}\| \widetilde \gamma '\|_{L^\infty} \leq 1 +
   \widetilde \kappa\big( \|\nabla \eta\|_{L^\infty}     \ell(T) + \ell(T) /\epsi_0 +1
   \big)=: \haz L <3/2\quad \text{a.e. in
     $(0,\ell_\gamma)$}.\label{eq:MM}
 \end{equation}
 A  similar argument entails that
 \begin{equation}\label{eq:slip}
   \| \widetilde
   \gamma'\|_{ L^\infty} \geq 1 - \widetilde \kappa \big( \|\nabla
   \eta\|_{L^\infty}      \ell(T)  +  \ell(T)  /\epsi_0 +1
   \big)> 1/2  \quad \text{a.e. in
     $(0,\ell_\gamma)$}.
   \end{equation}

 The lower bound \eqref{eq:elisa0} and the triangle inequality allow us to control
  \begin{align*} 
   &  \ove \kappa s \stackrel{\eqref{eq:elisa0}}{\leq} d (\gamma(s), \partial V(t)\setminus
    \partial \Omega)  \leq d(\widetilde\gamma(s), \partial V(t)\setminus
    \partial \Omega) + |\widetilde \gamma(s) - \gamma (s)|\nonumber\\
    &\quad =  d(\widetilde\gamma(s), \partial V(t)\setminus
    \partial \Omega) +  \big(\epsi_0 - d(\gamma(s) ,\partial \Omega)\big)^+ \frac{\widetilde \kappa}{\epsi_0} s
      \leq  d(\widetilde\gamma(s), \partial V(t)\setminus
    \partial \Omega) +  \widetilde \kappa  s 
  \end{align*}
  entailing that
  \begin{equation}
    \label{eq:john2}
    d (\widetilde \gamma(s), \partial V(t)\setminus  \partial \Omega)\geq (\ove
    \kappa -  \widetilde \kappa ) s> (2\widetilde
    \kappa -  \widetilde \kappa )  s = \widetilde \kappa s\quad \forall s \in [0,\ell_\gamma].
  \end{equation}

  Fix now $s\in [0,\ell_\gamma]$. If $d(\gamma(s), \partial \Omega) \geq \epsi_0$ we use the fact that $
  \widetilde \kappa < \epsi_0 /   \ell(T)  $ from \eqref{eq:tildek} to find $d(\widetilde
  \gamma(s), \partial \Omega) = d(\gamma(s),\partial \Omega)   \geq \epsi_0\geq    \widetilde \kappa    \ell(T)  \geq  
   \widetilde \kappa s$. On the contrary, if  $d(\gamma(s), \partial
   \Omega) < \epsi_0$ we have
   \begin{align*}
    &d(\widetilde \gamma (s), \partial \Omega) = d(\gamma(s),\partial \Omega) + 
      \big(\epsi_0 - d(\gamma(s), \partial  \Omega) \big)\frac{\widetilde \kappa}{\epsi_0}  s \nonumber\\
    &\quad = \left(1-\frac{\widetilde \kappa}{\epsi_0} s\right)
      d(\gamma(s),\partial \Omega) +  \widetilde \kappa  s\geq  \left(1-\frac{\widetilde \kappa}{\epsi_0}     \ell(T)  \right)
      d(\gamma(s),\partial \Omega) +  \widetilde \kappa
      s\stackrel{\widetilde \kappa<\epsi_0/    \ell(T)   }{\geq}  \widetilde \kappa s.
   \end{align*}
  In both cases, we hence have that 
  \begin{equation}
    \label{eq:john4}
    d(\widetilde \gamma (s), \partial \Omega) \geq  \widetilde \kappa  s \quad
    \forall s \in [0,\ell_\gamma].
  \end{equation}

  Putting together \eqref{eq:john2} and \eqref{eq:john4} we conclude that
  \begin{align}
    &d(\widetilde\gamma(s), \partial V(t))= \min\big\{d(\widetilde\gamma(s), \partial
  V (t) \setminus \partial \Omega), d(\widetilde\gamma(s),  \partial
      V(t)\cap \partial \Omega)\big
      \}\geq \widetilde
      \kappa s \quad \forall s \in [0,\ell_\gamma]. \label{eq:john5}
  \end{align}
  In particular, we have that $\widetilde \gamma (\cdot) \in V(t)
 $ 
 (recall that $\widetilde \gamma(0)=x \in V(t)$).  
\medskip

\noindent  {\it Step 3: Arc-length reparametrization of $\widetilde \gamma$.}

 The last step of the proof consists in reparametrizing $\widetilde \gamma$
 by arc length: One defines $\haz \gamma :[0,\ell_{\haz \gamma}]\to \ove
 \Omega$ as $\haz\gamma(\tau) = \widetilde \gamma(s(\tau))$ where $s: [0,\ell_{\haz
   \gamma}]\to [0,\ell_{ \gamma}]$ is the
 inverse function of $s\mapsto \tau(s):=\int_0^s |\widetilde \gamma'(r)|\, \d
 r$. Note that $s$ is Lipschitz continuous due to \eqref{eq:slip}. Using
 \eqref{eq:MM}, for all $\tau \in [0,\ell_{\haz \gamma}]$ we have that
 \begin{equation}
   s(\tau)= \int_0^\tau s'(t)\, \d t = \int_0^\tau \frac{\d
   t}{\tau'(s(t))} = \int_0^\tau \frac{\d
   t}{|\widetilde \gamma'(s(t))|} \stackrel{\eqref{eq:MM}}{\geq}\int_0^\tau \frac{\d
   t}{ \haz L} = \frac{\tau}{ \haz L}.\label{eq:esse}
 \end{equation}
 Hence, setting $\kappa
 : = \widetilde \kappa/ \haz L $ we obtain 
 $$d(\haz \gamma (\tau), \partial V(t)) =
 d(\widetilde\gamma(s(\tau)),\partial V(t))
 \stackrel{\eqref{eq:john5}}{\geq} \widetilde \kappa s(\tau) \stackrel{\eqref{eq:esse}}{\geq} \kappa
 \tau \quad \forall \tau \in [0,\ell_{\haz \gamma}].$$
 This implies that $V(t)$ is a John domain  with John center
 $x_0$ and John constant
 $$\kappa(T) = \frac{\widetilde \kappa}{1 + \widetilde \kappa
   (\| \nabla \eta \|_{L^\infty}  \ell(T)  +  \ell(T) /\epsi_0 +1)}>0.$$
 In particular, $\kappa(T)$ depends on $\kappa_0$, $\sigma_*$,
 $\sigma^*$, $\|\nabla \eta\|_{L^\infty}$,  $\ove
 \epsi$, $d(x_0,\partial \Omega)$,  $L$, $d(x_0,\partial V_0)$,
  and   $T$, and
 $\kappa(T) \to 0$ as $T \to \infty$.  

 \medskip
 
\noindent  {\it Step 4: Lipschitz continuity of $t \mapsto V(t)$ with
   respect to the Hausdorff distance}.
 
Assume that $ 0<s<t\leq T$,  choose $x \in V(t)\setminus \overline{V(s)}$, and
let $\gamma_x\in  {\mathcal C}_x $ be optimal for $x$ with $\gamma_x(\cdot) \in \ove
\Omega$. For $ \iota  >0$ small, fix
$\tau\in(0,1)$ with $ 0 <  s- \iota   <v(\gamma_x(\tau))<s$, so that $\gamma_x(\tau) \in V(s)$. Such $\tau$
exists, as $ r  \in [0,1]\mapsto v(\gamma_x( r ))$ is continuous,
 $v(\gamma_x(0))=0$, and $v(\gamma_x(1))=v(x)>s >0 $.  We have that
\begin{align*}
 & \sigma_* d(x,V(s)) \leq \sigma_* |x-\gamma_x(\tau)|\leq
   \sigma_*d_{\ove \Omega}(x,\gamma_x(\tau)) \leq \sigma_*\int_\tau^1
   |\gamma_x'(r)|\, \d r \\
  &\hspace{-2mm}\quad \stackrel{\eqref{eq:sigma}}{\leq}
   \int_\tau^1\sigma(\gamma_x(r),\gamma_x'(r))\, \d r
  =
    \int_0^1\sigma(\gamma_x(r),\gamma_x'(r))\, \d r-
     \int_0^\tau\sigma(\gamma_x(r),\gamma_x'(r))\, \d
    r \\[1mm]
  &\quad  \leq  v(x) - v(\gamma_x(\tau)) < t - s +  \iota .
\end{align*}  
By taking the infimum for $ \iota  >0$ and then the supremum on $x\in
V(t)$ we obtain \eqref{eq:lipschitz}.

\medskip

\noindent {\it Step 5: H\"older  continuity of $t \mapsto V(t)$ with
   respect to the $L^1$ distance}.
 
 Let $0<s<t\leq T$ and recall that $V( s )$ is a John domain with respect to
$x_0$ with John constant $\kappa(T)$. Owing to
\cite[Cor.~6.4]{Martio2}, see also \cite[Cor.~2]{Smith}, 
the Minkowski (or box-counting) dimension of $\partial V( s )$ is
smaller than $n$. More precisely,  there exists $r_0>0$ such that 
the minimal number $N(r)$ of balls
with radius $ r   \in (0,r_0) $  needed to
cover $\partial V( s)$ can be bounded as
\begin{equation}\label{eq:Martio}
  N(r)\leq c r ^{\mu-n}
  \end{equation}
  where  $r_0>0$,  $\mu\in (0,1)$ and $c>0$ depend on  $\ell(T)$
  from \eqref{eq:lmax}  
and $\kappa(T)$ only.

Fix $r:=2(t-s)/ \sigma_*$.  If $r \leq r_0$ consider a finite
covering $(B_r(x_j))_{j=1}^{N(r)}$ of $\partial V(s)$ with balls of
radius $r$ and define $ F  = \cup_{j=1}^{N(r)} B_{3r/2}(x_j)$.  The bound   
\eqref{eq:lipschitz} guarantees that 
$d_H(V(t),V(s))\leq (t-s)/\sigma_*=r/2$ so that 
$V(t)\setminus V(s) \subset  F $. Using \eqref{eq:Martio} we
conclude that 
$$|V(t)\setminus V(s) | \leq | F |\leq  N(r)|B_1|  \frac{3^n}{2^n}
 r^n
  \leq c |B_1|\frac{3^n}{2^n}\left(\frac{2(t-s)}{\sigma_*}
  \right)^{\mu} =  \frac{3^n c |B_1| }{2^{n-\mu}\sigma_*^\mu} (t-s)^{\mu} $$
  where $\mu\in (0,1)$ and $c$ depend on   $\ell(T)$   and the John constant $\kappa(T)$ of $V(t)$
only.

 If on the contrary $r\geq r_0$ we have that $1\leq 2(t-s)/
\sigma_*r_0 $, hence
\begin{align*}
  |V(t) \setminus V(s)| \leq |V(T)| \leq \left({\rm diam}(V_0) +
  {2T}/{\sigma_*}  \right)^n \leq \frac{2^\mu \left({\rm diam}(V_0) +
  {2T}/{\sigma_*}  \right)^n }{\sigma_*^\mu r_0^\mu}(t-s)^\mu.
\end{align*}
By letting $\mu(T)= \mu$ and
  $$c_{ H}(T) = \max\left\{  \frac{3^n c |B_1| }{2^{n-\mu}\sigma_*^\mu}, \frac{2^\mu\big({\rm diam}(V_0) +
     {2T}/{\sigma_*} \big)^n }{\sigma_*^\mu r_0^\mu}\right\}$$
we obtain \eqref{eq:hoelder}.
\end{proof}

\section{Free boundary problem: Proof of
  Theorem~\ref{thm:main}}\label{sec:elastic}

 We start by recording a uniform Poincar\'e estimate.
 
\begin{proposition}[Uniform Poincar\'e
  inequality]\label{prop:poincare}
  Assume \eqref{eq:H05}--\eqref{eq:H07}. For all
  $ \underline{\kappa} \in (0,\kappa_0)$ and $R>0$ define
\begin{align*}
&\notag\mathcal{V}:= \big\{ V  \subset \Rz^n \ | \ V_0 \subset      
                V \subset  \Omega \cap B_R\ \text{and $V$ is open
                and John regular with
                               respect to $x_0
                               $,}\\
                         & \qquad \qquad \text{with John constant
                           $\kappa\in [\underline{\kappa},\kappa_0]$}\}.
\end{align*}
Then, there exists a constant $
c_{\mathcal{V}}>0$ depending on $\underline{\kappa}$, $d(x_0,\partial
V_0)$,  $\Omega$, $n$,   and $R $ 
such that 
\begin{equation}
\label{eq:unif-korn}
\|  u \|_{H^1_\Gamma(V)} \leq c_{\mathcal{V}}\|\nabla u\|_{L^2(V;\Rz^n)}  \quad \forall V \in
 \mathcal{V}, \ \forall u \in H^1_\Gamma(V). 
\end{equation}
\end{proposition}

\begin{proof} Let $\haz x \in \Gamma $ and $\haz\epsi>0$ be such that
  $d(\haz x,\partial
  V_0\setminus \partial \Omega) >2\haz \epsi$   and $(B_{\haz
    \epsi}(\haz x) \cap \partial \Omega) \subset \Gamma$.  Take any $V \in
  {\mathcal V}$ and let $\kappa$ be its John constant. We claim that
  the enlarged set
  \begin{equation}
    \label{eq:claim}
   V_{\haz \epsi}:= V\cup B_{\haz \epsi}(\haz x) \quad \text{is John regular w.r.t. $x_0$ with
      John constant $\kappa/3$}.
  \end{equation}
  As soon as \eqref{eq:claim} is proved, one can apply 
  \cite[Lemma~3.1]{bojarski} and
\cite[Theorem~5.1]{bojarski} on the John domain $V_{\haz \epsi}$ to
functions in $  u
\in H^1(V_{\haz \epsi}) $ with $u=0$ on $V_{\haz
  \epsi}\setminus \ove{V}$ to obtain the  assertion.   In particular,
$u=0$ on $B_{\haz \epsi}(\haz x)\setminus \ove \Omega$ which has a
positive  measure as $\Omega$ is uniformly $C^{1,1}$, hence the
external-ball condition holds.  

We are hence left with proving claim \eqref{eq:claim}. This follows from a
simplified version of Step  1 in the proof of Proposition \ref{prop:john}. As $V$ is John regular,
one can find an arc-length parametrized curve
$\haz\gamma:[0,\haz\ell] \to \ove{V}$ with $\haz\gamma(0)=\haz x$,
$\haz \gamma(\haz\ell)=x_0$, and $d(\haz\gamma(s),\partial V)\geq \kappa s $ for all
$s \in [0, \haz\ell]$. Let now  $x \in V_{\haz \epsi}\setminus
\ove{V}$   (the proof for points in $\ove V$ being trivial)  and define the  arc-length parametrized curve $ \gamma:[0,\haz \ell+\ell]\to
V_{\haz \epsi}$ by setting $ \ell:=|\haz x - x|<\haz \epsi$ and 
$$
\gamma(s) =
\left\{
  \begin{array}{ll}
    x+s(\haz x -x)/\ell\quad &\text{for} \ s \in [0,
                                              \ell],\\[1mm]
    \haz \gamma(s- \ell)&\text{for} \ s \in (\ell, \haz \ell+
                               \ell] .
  \end{array}
  \right.
$$
Note that $\gamma(0)=x$ and $\gamma(\haz \ell+
\ell) = \haz \gamma(\haz \ell)=x_0$. One readily
 checks  that
\begin{equation}\label{eq:claim1}d(\gamma(s),\partial V_{\haz
    \epsi}) \geq s\quad \forall s \in [0,\ell].
\end{equation}
If $s\in
[\ell,3\ell/2]$ we compute that
$$|\haz x - \gamma(s) | = |\haz x- \haz \gamma(s-\ell)| = |\haz
\gamma(0) - \haz \gamma(s-\ell)| \leq s-\ell\leq \frac{\ell}{2} .$$
The triangle inequality hence implies that
\begin{align*}
  &d(\gamma(s), \partial V_{\haz \epsi}) \geq d(\haz x,  \partial
    V_{\haz \epsi}) - |\haz x - \gamma(s)| \nonumber\\
  &\quad  \geq  \haz \epsi - |\haz x - \gamma(s)| > \ell- \frac{\ell}{2}=
  \frac{\ell}{2} \geq \frac{s}{3}\quad \forall s \in [\ell,3\ell/2]. 
\end{align*}
 Using $V \subset V_{\haz \epsi}$  we deduce that
\begin{align}
  d(\gamma(s), \partial V_{\haz \epsi}) \geq  d(\gamma(s), \partial
  V ) = d(\haz \gamma (s -\ell),\partial V) \geq \kappa
  (s-\ell)\geq \frac{\kappa}{3}s \quad \forall s\in [3\ell/2,\haz \ell+\ell].\label{eq:claim3}
\end{align}
By collecting the information in \eqref{eq:claim1}--\eqref{eq:claim3}
and recalling that $\kappa \leq \kappa_0\leq 1$,  claim
\eqref{eq:claim} follows.
\end{proof}

\begin{proof}[Proof of Theorem~\ref{thm:main}]
 We  use a fixed-point argument,   which   we divide into
  steps. At first, we  consider the elliptic 
  problem \eqref{eq:0111} for   given growing sets  (Step 1).  Then, we
  solve the Hamilton--Jacobi subproblem
  \eqref{eq:0222} for a given  $Ku(x,v(x))$  (Step 2). This
  allows us to define  a map $S$  
(Step 3)   which is proved to admit a fixed point  (Step
   4). Finally,  we prove that  
  $v$ is a constrained viscosity solution  (Step  5). 
  \medskip

  \noindent {\it Step 1:  Solution to the elliptic   problem, given the growth.} 
   Define the set
  \begin{align}
    &M:= \{  v \in C(\ove \Omega)\ |  \ \text{$\forall T>0, \ t\in [0,T]$ the sublevel set}\nonumber\\
    &\quad \text{$ V(t)=\{x\in  \Omega \ | \
       v(x)<t\}$ contains $V_0$ and is a  John domain}\nonumber\\
    &\quad \text{with respect to $x_0 \in V_0$
      with John constant $\kappa(T) $}\}.\label{eq:M}
  \end{align}
  Assume to be given $ v \in M$. Owing to the standard Lax--Milgram
  Lemma, for all fixed $t >0$ one can
  find a unique  $u(\cdot,t) \in H^1_\Gamma(V(t))$ solving
  \begin{equation}
    \int_{V(t)}\nabla u(x,t)\cdot \nabla w (x)\, \d x =
  \int_{V(t)}w(x)\, \d x \quad \forall w \in
  H^1_\Gamma(V(t)).\label{eq:pratico}
  \end{equation}
  Moreover, one can prove that the ensuing $u:Q_{ v}=\cup_{t>0}
  V(t)\times \{t\}\to \Rz $ is measurable.  In fact, one can check
  that the map $t \mapsto \ove u(\cdot ,t)$ is weakly left-continuous
  in $L^2(\Rz^n)$ at all times  as $t \mapsto V(t)$ increases.   By choosing $w=u(\cdot,t)$
in \eqref{eq:pratico}  and using the uniform Poincar\'e inequality
\eqref{eq:unif-korn}  with $R=R(t) = {\rm diam} V_0 +2t/\sigma_*$  one gets that 
\begin{align}
  \label{eq:estimateu0}
  \|\ove u(\cdot,t)\|_{L^2(\Rz^n)}&\leq c_{ 2}( t ) \quad
   \forall t >0,\\
  \|u(\cdot,t)\|_{H^1_\Gamma(V(t))}&\leq c_{ 2}( t ) \quad
   \forall t >0. \label{eq:estimateu}
\end{align}
 for some increasing function $c_{ 2}:\Rz_+\to\Rz_+$  depending
 just on data.  
\medskip

\noindent {\it Step 2: Solution to the Hamilton--Jacobi problem, given
  $K u(x, v(x))$.}   Let  $\tilde v\in
C(\ove \Omega)$  be given  and $ \tilde \zz \in
L^\infty_{\rm loc}( [0,\infty) ;L^2(\Rz^n))$, so that $K \tilde \zz (\cdot, \tilde
v(\cdot)) \in C( \ove \Omega )$. Defining $$\tilde H(x,p):= H(x,
K \tilde \zz (x, \tilde
v(x)) ,p)$$ for all $ x \in \ove \Omega$ and $p\in \Rz^n$, we have that
$\tilde H$ fulfills the assumptions
\eqref{eq:H01}--\eqref{eq:H04}. In fact, \eqref{eq:H01}, \eqref{eq:H02}, and
\eqref{eq:H04} follow directly from the continuity of $K \tilde \zz (\cdot,
\tilde v(\cdot))$ and \eqref{eq:nH01u}, \eqref{eq:nH02u}, and
\eqref{eq:nH04u}, respectively. The remaining
\eqref{eq:H03}--\eqref{eq:H03b} are readily checked via 
\eqref{eq:H03u}--\eqref{eq:H03ub}. This allows us to apply 
Proposition~\ref{prop:representation} and find a
viscosity solution to
\begin{align*}
  &H(x, K \tilde \zz (x,\tilde v(x)),\nabla v )=0 \quad \text{in} \ \Omega\setminus
  \ove{V_0}, \\
  &v=0 \quad \text{in} \ V_0
\end{align*}
given by
\begin{equation*}v(x)= \min \Bigg\{\int_0^{1} 
  \sigma\big(\gamma,K \tilde \zz (\gamma,\tilde v (\gamma)),\gamma'\big) \
  {\d s} \ \Big| \   \gamma \in  {\mathcal C}_x , \ \gamma (\cdot)\in
  \ove {\Omega}\Bigg\} \quad \forall  x \in \ove{\Omega}.
\end{equation*}
Note that Propositions \ref{prop:estimates}--\ref{prop:john}
  guarantee that $v \in M$, see \eqref{eq:M}. 
\medskip

\noindent {\it Step 3: Fixed-point map.} Endow $E= L^\infty_{\rm
  loc}(\Rz_+;L^2(\Rz^n)) \times C(\ove \Omega)$ with the topology $\tau$ of local
weak* convergence $\times$ local uniform convergence, which makes 
$E$  a  locally convex Hausdorff  
space. Define $$W:=\{( \zz ,v) \in E \  | \, \| \zz  (\cdot,t)\|_{L^2(\Rz^n)} \leq c_2(t) \ \text{for a.e.} \  t>0, \ |\nabla v| \leq
\sigma^*L \  \text{a.e.},  \ v =0 \ \text{on} \ \ove{V_0}\}$$
which is not empty, convex, and compact with respect to
$\tau$.  Note that $\tau$ is metrizable on~$W$. 

Define the set-valued map $S:W \to 2^E$ by setting $ ( \zz , v)\in S( \tilde \zz , \tilde v) $ iff
\begin{align}
  &v(x) = \min \Bigg\{\int_0^{1} 
  \sigma\big(\gamma,K \tilde \zz (\gamma, \tilde v (\gamma)),\gamma'\big) \
  {\d s} \ \Big| \   \gamma \in  {\mathcal C}_x , \ \gamma (\cdot)\in
  \ove {\Omega}\Bigg\} \quad \forall x \in \ove
    \Omega, \label{eq:S5}\\
    &H(x, K \tilde \zz (x,\tilde v(x)),\nabla v )=0 \quad \text{in} \ \Omega\setminus
      \ove{V_0}, \label{eq:1js0}\\[2mm]
  & \int_{V(t)} \nabla  \zz (\cdot,t)\cdot \nabla w \, \d x =
    \int_{V(t)}w \, \d x \quad \forall w\in H^1_\Gamma(\Omega), \
    \text{for  a.e.} \ t>0 \nonumber\\
  & \qquad \qquad \text{where} \ \ V(t) =
     \{x\in \Omega \ | \ v(x)<t\} \quad \forall t>0, \label{eq:S00}\\[2mm]
  &\|  \zz (\cdot,t) \|_{L^2(\Rz^n)} \leq c_2(t) \quad \text{for a.e.} \  t>0,\label{eq:S1}\\[2.5mm]
  &\| \zz (\cdot,t)\|_{H^1(V(t))} \leq c_2(t) \quad \text{for a.e.} \ 
    t>0, \label{eq:S4}\\[2mm]
  & \zz (\cdot, t) \in H^1_\Gamma(V(t)) \quad \quad \text{for a.e.} \  t>0,\label{eq:1js1}\\[2mm]
  & \zz =0 \ \text{in}  \ (\ove {Q_{ v}})^c .\label{eq:S2} 
\end{align}

At first, let us observe that $S$  takes   nonempty values in $W$. In fact, Step 2
ensures that $v$ defined by \eqref{eq:S5} solves \eqref{eq:1js0} in
the viscosity sense, $v =0 $ in $V_0$, and $v\in M$. Hence, owing to
Step 1 one can find $ \zz  $ fulfilling
\eqref{eq:S00}--\eqref{eq:S2}.  In particular, as $v \in M$ and
\eqref{eq:estimateu0} holds, we have that $( \zz ,v)
\in W$, as well.
Note that  $ \zz $ from Step 1 solves the
stronger \eqref{eq:pratico}, where tests are taken in $
H^1_\Gamma(V(t))$. Let us stress again that \eqref{eq:pratico}
and \eqref{eq:S00} are
not equivalent if $V(t)$ is not an extension domain. In particular,
 (restriction to $V(t)$ of) 
the solution to \eqref{eq:pratico} is unique, while  
\eqref{eq:S00}--\eqref{eq:S2} imply no uniqueness, in general.

It is a standard matter to prove that $S$ has convex values. Indeed,
\eqref{eq:S5} uniquely defines $v$. On the other hand, relations
\eqref{eq:S00}, \eqref{eq:1js1}, and \eqref{eq:S2} are  affine,  while
conditions \eqref{eq:S1}--\eqref{eq:S4} are convex.

In Step 4 below, we provide a proof of the following lemma.

\begin{lemma}[Closed graph]\label{lemma:closed}  $\{((\tilde \zz,\tilde v) , (\zz,v)) \in W \times W \  | \
(\zz,v)\in S(\tilde \zz,\tilde v)\}$ is closed. 
\end{lemma}

Owing to this lemma, we are in the position  to apply  the Kakutani--Fan--Glicksberg
Theorem \cite[Thm. 1.11, p.~8]{roubicek} which ensures the existence
of a fixed point  $(\zz,v)$ for $S$. By letting $u$ be the restriction of
$\zz$ to $Q_v$, we have that $(u,v) $ solves \eqref{eq:01} and
\eqref{eq:0111}--\eqref{eq:0222} and fulfills \eqref{eq:S4}--\eqref{eq:1js1}. 
This
proves the existence part of Theorem \ref{thm:main}.
\medskip

\noindent{\it Step 4: Proof of Lemma \ref{lemma:closed}.} 
In order to check the  closedness  of the graph of $S$,
recall that the topology $\tau$ is metrizable  on $W$  and take
$( \tilde \zz _j,\tilde v_j)\stackrel{\tau}{\to}(\tilde
{ \zz }, \tilde v)$ and $ ( \zz _j,v_j)\in S( \tilde \zz _j,\tilde v_j)$ with
$( \zz _j,v_j) \stackrel{\tau}{\to}( \zz ,v)$.  Lemma
\ref{lemma:closed} follows by  proving   that $( \zz ,v)\in S(\tilde
{ \zz },
\tilde v)$, i.e., that \eqref{eq:S5}--\eqref{eq:S2} hold.

Ad \eqref{eq:S5}--\eqref{eq:1js0}:  As $ \tilde \zz _j \to  \tilde \zz $
weakly* in $L^\infty_{\rm loc}( [0,\infty) ;L^2(\Rz^n))$, by applying the mollification we readily check that 
\begin{align}
  &K  \tilde \zz _j\to  K { \tilde \zz } \quad \text{locally
    uniformly in} \  \Rz^n \times \Rz_+.
    \label{eq:pj900}
\end{align}
In addition, using the fact that
$$\partial_t K  \tilde \zz _j (x,t)= k(0)
\int_{\Rz^n}\phi(x-y) \tilde \zz _j (y,t)\, \d
y + \int_0^t\!\!\int_{\Rz^n} k'(t-s)\phi(x-y) \tilde \zz _j(y,s)\, \d y \, \d s$$
we can prove that 
$K\tilde \zz_j$ are locally Lipschitz
continuous in time, namely, 
\begin{align*}
 &\| K \tilde \zz _j(\cdot,t_1) -K \tilde \zz _j(\cdot  ,t_2)
  \|_{L^\infty(\Rz^n)}  \leq \|
  \partial_t   K  \tilde \zz _j\|_{L^\infty(\Rz^n \times (
   0, T  ))} |t_1 - t_2| 
  \nonumber \\ 
  & \quad  \leq \big(  |k(0)|+\| k'\|_{L^1(
   0, T  )}\big)
    \|\phi\|_{L^2(\Rz^n)}   \| \tilde \zz _j
    \|_{L^\infty(
   0, T  ;L^2(\Rz^n))} |t_1 - t_2|\nonumber\\
 &\quad \leq \big(  |k(0)|+\| k'\|_{L^1(
   0, T  )}\big)
    \|\phi\|_{L^2(\Rz^n)}   c_{ 2}(T)|t_1 - t_2|
\end{align*}
for all $0<t_1 <t_2<T$.  Hence, for   all $x \in \Omega$ we can control
\begin{align*}
  & |K  \tilde \zz _j (x,\tilde v_{j}(x))
    - K \tilde \zz (x,\tilde v(x))|\nonumber \\
  &\quad \leq  |K  \tilde \zz _j
  (x,\tilde v_{j}(x)) - K  \tilde \zz _j (x,\tilde v(x))|  +  |K
     \tilde \zz _j (x,\tilde v(x)) -
    K \tilde \zz (x,\tilde v(x))| \nonumber\\
  &\quad\leq \big( |k(0)|+\| k'\|_{L^1(
   0, T  )} \big)
    \|\phi\|_{L^2(\Rz^n)}   c_{ 2}(T)\, |\tilde v(x) -\tilde v_j(x)| \nonumber\\[2mm]
  &\qquad + |K \tilde \zz _j (x,\tilde v(x))- K \tilde \zz (x,\tilde v(x))|.
\end{align*}
 As $\tilde v_j \to \tilde v$ locally uniformly and
convergence \eqref{eq:pj900} holds, this proves that  
$K \tilde \zz _j(\cdot,\tilde v_j(\cdot)) \to K\tilde \zz(\cdot,\tilde v(\cdot)) $
locally uniformly.  By standard arguments (see, e.g., Step 4 in
the proof of Proposition~\ref{prop:representation}) this  implies that $v$ is a viscosity solution
to \eqref{eq:1js0}. 

In order to check the representation formula \eqref{eq:S5}, set $x \in \ove
\Omega$. For any $\gamma\in  {\mathcal C}_x $ with $\gamma(\cdot )\in \ove \Omega$, using
the local uniform convergence of $ K \tilde \zz _j(\cdot,\tilde
v_j(\cdot))$ proved above we obtain that 
\begin{equation*}
  K \tilde \zz _{j}(\gamma,\tilde v_{ j}(\gamma)) \to
  K \tilde \zz (\gamma,\tilde v(\gamma))\quad \text{uniformly in} \ [0,1].
\end{equation*}
By the continuity
of $\sigma$  (which holds as in Lemma \ref{lemma:sigma} also for
$\sigma$ additionally depending on $u$)  and by dominated
convergence we hence have that
\begin{align}
 v(x)= \lim_{j\to \infty}v_j(x)\leq \lim_{j\to
  \infty}\int_0^1
  \sigma(\gamma,K   \tilde \zz _{j}(\gamma, \tilde v_{j}  (\gamma)),\gamma')\, \d
  s=\int_0^1
  \sigma(\gamma,K  \tilde \zz  (\gamma,\widetilde v(\gamma)),\gamma')\, \d
  s.\label{eq:pj110}
\end{align}
To prove \eqref{eq:S5}, we need to find an optimal curve realizing the equality in \eqref{eq:pj110}.
We argue as in the
proof of Proposition~\ref{prop:representation}. Let
$\gamma_{xj} \in  {\mathcal C}_x $ with $\gamma_{xj}(\cdot)\in \ove \Omega$ be optimal curves for $x$
under $\sigma(\cdot,K \tilde \zz _{j}(\cdot,  \tilde v_{j}
(\cdot)),\cdot)$. Letting  $\ell_j$ be the length of $\gamma_{xj}$, define
 $s_j:[0,1]\to [0,1+\ell_j]$   to  be $s_j(t):=t+ \int_0^t
 |\gamma_{xj}'( r)|\, \d  r$ for $t
\in [0,1]$, 
denote by $t_j: [0, 1+\ell_j] \to [0,1] $ its inverse,
and define $\haz \gamma_{xj}( r):=
\gamma_{xj}(t_j((1+\ell_j)  r))$ for $ r \in
[0,1]$. Then, $\haz \gamma_{xj}\in  {\mathcal C}_x $  is optimal  and $|\haz
\gamma_{xj}'| \leq 1+ \ell_j $ a.e. Moreover, $\haz
\gamma_{xj}(1)=x$ and $\tilde v_j(x)\to
\tilde v(x)$, so that the sequence $(\haz \gamma_{xj})_j$ is uniformly
bounded in $W^{1,\infty}(0,1;\Rz^n)$. By extracting without relabeling
we have  that  $\haz \gamma_{xj}\to \gamma$
uniformly in $[0,1]$ and weakly* in  $W^{1,\infty}(0,1;\Rz^n)$,  so
that $\gamma(\cdot) \in \ove \Omega$.   In
particular, $\gamma(1)=\lim_{j \to \infty}\haz \gamma_{xj}(1)=x$
and  $ \gamma(0) \in
\ove{V_0}$ as 
$\ove{V_0}\ni \haz \gamma_{xj}(0) \to \gamma(0)  $ and
$\ove{V_0}$ is closed,  hence  $\gamma \in
 {\mathcal C}_x $. We hence conclude that
\begin{align*}
v(x)&=\lim_{j \to \infty}
                 v_{j}(x)=\lim_{j \to \infty}\int_0^{1}\sigma
      (\haz\gamma_{xj}(s), K \tilde \zz _{j}(\haz\gamma_{xj}(s),  \tilde v_{j}  (\haz\gamma_{xj}(s))), \haz\gamma_{xj} '(s))\,
           \d s  \\
  &\geq    \int_0^{1}\sigma (\gamma(s),  K\tilde
    { \zz }
    (\gamma (s),\tilde v  (\gamma (s))), \gamma  '(s))\,
           \d s  
\end{align*}
by \cite[Thm.~3.23, p.~96]{Dacorogna}, as $\sigma$ is
nonnegative, continuous, and convex in $\gamma'$.  This proves~\eqref{eq:S5}. 

Before moving on, let us remark that Proposition \ref{prop:john} holds
for $V(t)=\{x \in
\Omega \ | \ v(x) <t\}$ and $V_j(t)=\{x \in
\Omega \ | \ v_j(x) <t\}$. In particular, $V(t)$ and $V_j(t)$ are John
domains with respect to the same $\kappa(\cdot)$.

Ad \eqref{eq:S00}:  As $(\zz_j)_j $ fulfill \eqref{eq:S4}, 
setting $\hat \zz_j$ to be the restriction of $\zz_j$ to $Q_{v_j}$ and
extracting a not relabeled subsequence  we have
  \begin{align} 
    &\ove{\nabla  { \hat \zz }_j} \to  S  \quad \text{weakly* in} \
      L^\infty_{\rm loc}( [0,\infty) ;L^2(\Rz^n;\Rz^n))\label{eq:pj50} 
  \end{align}
  where $\ove{\nabla  { \hat \zz }_j}(\cdot, t) $ denotes the trivial extension
  of ${\nabla  { \hat \zz }_j} (\cdot, t) $ from $V(t)$ to $\Rz^n$ for a.e. $t>0$.
  As $ v_j \to  v$ locally uniformly,
  for all 
$(x,t)\in  \Omega  \times \Rz_+$ with $ v(x)<t$ ($ v(x)>t$,
respectively), one has that $ v_j(x) <t$ ($ v_j(x)>t$) for $j$
large enough. In particular, the characteristic functions  $(x,t)\in \ove \Omega \times \Rz_+\mapsto I_j(x,t):=
1_{V_j(t)}(x)$  converge pointwise to $(x,t)\in \ove \Omega \times \Rz_+\to
I(x,t)=1_{V(t)}(x)$ at all points $(x,t)$ with $ v(x)\not = t$. On
the other hand, the set $\{(x,t)\in \ove \Omega\times \Rz_+\:| \: v(x)=t\}$
is negligible, so that the above convergence holds almost
everywhere. We conclude that
\begin{equation}
  \label{eq:pj700}
  I_j \to I \quad \text{strongly in} \ L^p_{\rm loc}( \Omega\times \Rz_+) \quad
  \forall p\in[1,\infty).
\end{equation}

Fix now an arbitrary $(x,t)\in  Q_{ v} $ and let $\delta>0$ be small, in such a
way that $A:= B_\delta(x)\times  (t-\delta,t+\delta)\subset \subset
 Q_v $.
Owing to the local uniform convergence of $ v_j$
 one has that $A\subset Q_{  v_j}$ for
$j$ large enough. By restricting to $A$ we hence have that
$ { \zz }_j\to  { \zz }$ weakly* in
$L^\infty(t-\delta,t+\delta;H^1(B_\delta(x)))$, which entails that
$  S =\nabla  \zz $ in $A$. By exhausting $Q_v$, this proves that 
\begin{equation} \label{eq:G0}
    S =\nabla  \zz  \quad \text{a.e. in} \ Q_v. 
 \end{equation}
  Similarly we prove that $S=0$ a.e. in $(\Rz^n \times
 (0,\infty))\setminus \ove{Q_v}$. 
 
  Let $D\subset H^1_\Gamma(\Omega)$ be dense and countable. For
all $w \in D$, let $J_w\subset \Rz_+$ be  the null set of
non-Lebesgue  points of
\begin{equation}
  t \in \Rz_+ \mapsto \int_{V(t)}\nabla  \zz (x,t)\cdot \nabla
  w(x)\,\d x=  \int_{\Rz^n}\ove{\nabla  \hat \zz }(x,t)\cdot \nabla w(x)\,
I(x,t)\,\d x\label{eq:lebesgue}
\end{equation}
 where $\hat \zz$ is the restriction of $\zz$ to $Q_v$ and
$\ove{\nabla  \hat \zz }$ is the trivial extension of $\nabla 
\hat \zz$ to $\Rz^n$.  
Set $ t \in \Rz_+\setminus J$ for $J:=\cup_{w\in
  D} J_w$, choose $w \in D$ in \eqref{eq:S00} written for $ \zz _j$,  and  
integrate over
$(t-\delta,t)$ for some small $ \delta \in (0,t)$  getting 
\begin{align*}
  &\frac{1}{\delta}\int_{t-\delta}^t \! \int_{\Rz^n} \ove{\nabla
   \hat \zz _j}(x,s)\cdot \nabla w(x)\, I_j (x,s)\, \d x \, \d s = \frac{1}{\delta}\int_{t-\delta}^t \! \int_{V_j(s)}{\nabla
   \zz _j}(x,s)\cdot \nabla w(x)\, \d x \, \d s \\
  &\quad =
  \frac{1}{\delta}\int_{t-\delta}^t \!\int_{V_j(s)} w(x)\, \d x
  \, \d s = \frac{1}{\delta}\int_{t-\delta}^t \!\int_{\Rz^n}  w(x)\, I_j (x,s)\, \d x \, \d s.
\end{align*}
Convergences
\eqref{eq:pj50} and \eqref{eq:pj700}  and the identification
\eqref{eq:G0}  allow us to pass to the limit  as $j \to \infty$  and get
\begin{align}
  &\frac{1}{\delta}\int_{t-\delta}^t \!\int_{V(s)} \nabla
     \zz  (x,s)\cdot \nabla w(x)\, \d x \, \d s
    =\frac{1}{\delta}\int_{t-\delta}^t \! \int_{\Rz^n} \ove{\nabla
    \hat \zz} (x,s)\cdot
    \nabla w(x)\, I(x,s)\, \d x \, \d s  \nonumber\\
  &\quad =\frac{1}{\delta}\int_{t-\delta}^t \!\int_{\Rz^n} w(x)\, I(x,s)\, \d x
  \, \d s= 
  \frac{1}{\delta}\int_{t-\delta}^t \!\int_{V(s)} w(x)\, \d x
  \, \d s.\label{eq:pj80}
  \end{align}
In order to pass to the limit as $\delta \to 0$, one can handle the
above right-hand side as follows
\begin{align*}
  &\left|\frac{1}{\delta}\int_{t-\delta}^t \!\int_{V(s)} w(x)\,
    \d x\,  \d s   -   \int_{V(t)} w(x)\, \d x
    \right| = \left| \frac{1}{\delta}\int_{t-\delta}^t
    \!\int_{V(t) \setminus V(s)}w(x)\, \d x\, \d s\right|\\
  &\quad \leq 
  \frac{1}{\delta}\int_{t-\delta}^t
    \!\int_{V(t)\setminus V(s)} |w(x)|\, \d x \, \d s\leq 
    |V(t)\setminus V(t-\delta)|^{1/2}\|w\|_{L^2(\Omega)}\\
  &\hspace{2mm}\stackrel{\eqref{eq:hoelder}}{\leq}
    (c_H(t)\delta ^{\mu(t)})^{1/2}\|w\|_{L^2(\Omega)} \to 0
\end{align*}
where we have used Proposition \ref{prop:john}.
As $t$ is a Lebesgue point  from  \eqref{eq:lebesgue},
 passing to the limit as $\delta \to 0$ in \eqref{eq:pj80} we 
  get
  \begin{equation*}  \int_{V(t)}\nabla  \zz (x,t)\cdot
\nabla w(x) \,\d x =
   \int_{V(t)} w(x)\, \d x. 
 \end{equation*}
Relation \eqref{eq:S00} follows from the density of $D$  in
$H^1_\Gamma(\Omega)$   by recalling
that $|J|=0$.

Ad \eqref{eq:S1}--\eqref{eq:S4}: Condition \eqref{eq:S1} is 
immediate from the lower semicontinuity of the norm  with respect to the weak*
convergence in $L^\infty_{\rm loc}( [0,\infty)
;H^1(\Rz^n))$. For all $t, \, \delta>0$ we have that

\begin{align*}
  &\int_t^{t+\delta}\| z(\cdot, r)\|^2_{H^1(V(r))}\, \d r
  \stackrel{\eqref{eq:G0}}{=} \int_t^{t+\delta}\int_{\Rz^n}
  \left(|\ove z(x,r)|^2 + |S(x,r)|^2
  \right) \, \d x \, \d r \\
  &\quad \leq \liminf_{j\to \infty}\int_t^{t+\delta}\int_{\Rz^n}
  \left(|z_j(x,r)|^2 + |\ove{\nabla \haz z_j}(x,r)|^2
    \right) \, \d x \, \d r \\
  &\quad =  \liminf_{j\to
  \infty}\int_t^{t+\delta}\|z_j(\cdot,r)\|^2_{H^1(V(r))}\, \d r \leq \int_t^{t+\delta}c_2^2(r)\, \d r.
\end{align*}

This proves that $t \mapsto  \| \zz (\cdot,t)\|^2_{H^1(V(t))} 
\in L^\infty_{\rm loc}( [0,\infty) )$. Letting $t>0$ be a
Lebesgue point of $r\mapsto c_2(r)\in L^\infty_{\rm
  loc}(\Rz_+)$ and $\delta \to 0$, bound \eqref{eq:S4} follows.

Ad \eqref{eq:1js1}:  By restricting to the set $G\subset V_0$, bound 
\eqref{eq:S4} ensures that  $z_j \to z$ weakly* in $L^\infty_{\rm
  loc}([0,\infty);H^1(G))$.    In particular, $\zz(\cdot, t)=0$ on $\Gamma$ for
  a.e. $t>0$ and \eqref{eq:1js1} holds.

Ad \eqref{eq:S2}: Take $(x,t) \in (\ove{Q_{v}})^c$ and choose
$\delta>0$ small, so that $ A:=B_\delta(x)\times
(t-\delta,t +\delta) \subset \subset  (\ove{Q_{v}})^c$. As $v_j \to v$
locally uniformly, we have that $ A\subset  (\ove{Q_{v_j}})^c$ for
$j$ large enough. In particular, $ \zz _j=0$ on $ A$ for $j$
large enough. Passing to the limit as $j \to \infty$ this implies that
$ \zz =0$ on $ A$, proving \eqref{eq:S2}.

\medskip
\noindent {\it Step 5:  Boundary behavior}. Once the representation
formula \eqref{eq:0222} is established, we can follow the argument of
Step 6 in the proof of Proposition~\ref{prop:representation}, by
additionally taking into account the dependence of $H$ on $Ku(\cdot,
v(\cdot))$.   Within  the contradiction
argument,   
using \eqref{eq:H03ub} 
we prove that there exists $\lambda\in (0,1)$ such that 
\begin{equation*}
  \gamma'(s)\cdot
\nabla \varphi(\gamma(s)) \leq \lambda \sigma
(\gamma(s),Ku(\gamma(s),v(\gamma(s))),\gamma'(s))\quad \text{for a.e.} \ s \in (t,1).
\end{equation*}
where $\varphi
\in C^1(\Rz^n)$ is such that $\varphi - v$ has a local maximum at
$x \in \partial \Omega \setminus \ove{V_0}$, $\gamma\in  {\mathcal C}_x $
realizes the minimum in \eqref{eq:0222}, $t\in (0,1)$, and $\gamma(s)
\in B_r(x)$ for $s \in [t,1]$  where  $r>0$ is fixed but arbitrary. A
straightforward adaptation of the chain
of inequalities \eqref{eq:can} yields
that $\varphi(\gamma(t)) - v(\gamma(t)) > \varphi(x)
- v(x)$, contradicting the fact that $x$ is a local maximum
of $\varphi - v$.

This concludes the proof  of Theorem \ref{thm:main}. 
\end{proof}

\section*{Acknowledgments and statements} 
This research was funded in
whole or in part by the Austrian Science Fund (FWF) projects
10.55776/F65,  10.55776/I5149, and 10.55776/PAT1408325. Part of this
work was developed
during the Thematic Program on {\it Free Boundary Problems} at the
{\it Erwin Schr\"odinger International Institute for Mathematics and
  Physics} in Vienna. The warm hospitality of the Institute is gratefully acknowledged. For
open-access purposes, the author has applied a CC BY public copyright
license to any author-accepted manuscript version arising from this
submission. A large language model (Fable 5.2, Nov. 2026) was used for proofreading.  The author assumes responsibility for
all content.

\end{document}